\documentclass[12pt,twoside,a4paper,draft]{amsart}
\usepackage{amssymb}
\usepackage{latexsym}
\usepackage{amsmath}
\usepackage{amsfonts}
\usepackage{dsfont}
\usepackage[OT4]{fontenc}
\usepackage{pst-node}
\usepackage{pst-plot}
\usepackage{pst-tree}

\newtheorem{tw}[equation]{Theorem}
\newtheorem{f}[equation]{Fact}
\newtheorem{w}[equation]{Corollary}
\newtheorem{lem}[equation]{Lemma}

\numberwithin{equation}{section}

\newenvironment{dowod}{\par\noindent {\bf Proof.}}
                      {\begin{flushright} \vspace*{-6mm} \mbox{$\square$} \end{flushright}\par}

\newcommand{\im}{\textnormal{Im}\,}
\newcommand{\h}{\mathcal H}
\newcommand{\F}{\mathcal F}

\newcommand{\R}{\mathbb R}
\newcommand{\C}{\mathbb C}
\newcommand{\N}{\mathbb N_0}
\newcommand{\la}{\lambda}

\newcommand{\de}{\delta}
\newcommand{\al}{\alpha}
\newcommand{\pd}{\sqrt{d}}
\newcommand{\gd}{\Gamma_d}
\newcommand{\ld}{\Lambda_d}

\newcommand{\bgd}{\partial\Gamma_d}
\newcommand{\dw}{d\omega}
\newcommand{\Ldwa}{L^2(\Omega, d\omega) }
\newcommand{\ldwa}{\ell^2(\Gamma_d) }
\newcommand{\ltwo}{\ell^2(\Lambda_d) }

\DeclareMathOperator{\obraz}{Im}

\DeclareMathOperator{\lin}{lin}
\DeclareMathOperator{\dist}{dist}
\DeclareMathOperator{\supp}{supp}

\newcommand{\1}[1]{\mathds{1}_{#1}}

\newcommand{\iloczyn}[2]{(\, {#1}\, , \, {#2}\, )}

\begin{document}

\title[Jacobi matrices on trees]
{Jacobi matrices on trees}
\author[A. M. Kazun and R. Szwarc]
{Agnieszka M. Kazun and Ryszard Szwarc}
\address{A. M. Kazun, Institute of Mathematics,
University of Wroc\l aw, pl.\ Grunwaldzki 2/4, 50-384 Wroc\l aw, Poland}
\email{agnieszka.kazun@gmail.com}
\address{R. Szwarc, Institute of Mathematics,
University of Wroc\l aw, pl.\ Grunwaldzki 2/4, 50-384 Wroc\l aw, Poland \newline and
\newline \indent Institute of Mathematics and Computer Science, University of Opole, ul. Oleska 48,
45-052 Opole, Poland}
\email{szwarc2@gmail.com}

\thanks{This research project has been partially supported by
European Commission Marie Curie Host Fellowship for the Transfer of Knowledge
``Harmonic Analysis, Nonlinear Analysis and Probability''
MTKD-CT-2004-013389.
Both the first and the second author were also supported by MNiSW Grant N201 3950 33 and
Grant N201 054 32/4285 respectively.
This work is based on the PhD thesis of the first author defended
in the Institute of Mathematics at the University of Wroc{\l}aw  in October 2008.}

\begin{flushright}
{\emph{
This work is dedicated \\
to Andrzej Hulanicki\\
our teacher and friend
}}
\end{flushright}

\hspace{14pt}

\maketitle

\thispagestyle{empty}

\begin{abstract}
Symmetric Jacobi matrices on one sided homogeneous trees are studied.
Essential selfadjointness of these matrices turns out to depend on the
structure of the tree. If a tree has one end and infinitely many origin points the matrix is always
essentially selfadjoint independently of the growth of its coefficients.
In case a tree has one origin and infinitely many ends, the essential
selfadjointness is equivalent to that of an ordinary Jacobi matrix
obtained by the restriction to the so called radial functions.
For nonselfadjoint matrices the defect spaces are described in terms of
the Poisson kernel associated with the boundary of the tree.
\end{abstract}

\hspace{7pt}

\begin{center}
{\textsc{INTRODUCTION}}
\end{center}

\hspace{7pt}

The classical moment problem consists in the following.

\hspace{7pt}

\noindent
\emph{
Given a sequence of real numbers $m_n.$ We look for a positive bounded measure $\mu$ on
the half-line $[0,\infty)$ or on the whole real line such that
\[
m_n\ = \ \int x^n d\mu (x)
\]
for $n=0,1,2,\ldots$.
}

\hspace{7pt}

Two main issues of the moment problem are existence and uniqueness of the measure $\mu$.
It is known that such a measure $\mu$ on the real line exists if and only if the numbers $m_n$ form
a  positive definite sequence.
A question of uniqueness of the measure $\mu$ is closely related to the selfadjointness of some operators.
The problem was intensively investigated
starting with the work of Thomas Jan Stieltjes (1894, \cite{stieltjes}, the case of the half-line)
and Hans Hamburger (1920, 1921, \cite{hamburger}, the case of the real line),
through that of Marcel Riesz (1921--23, \cite{riesz1}, \cite{riesz2}, \cite{riesz3}, a functional analysis approach),
Rolf Nevanlinna (1922, \cite{nevanlinna}, a complex function approach)
and Marshall H. Stone (1932, \cite{stone}, the Hilbert space methods),
until   recent results of Barry Simon (cf., e.g. \cite{simon}, 1998).

\hspace{7pt}

One of the key concepts that have arisen in the modern investigations is that of the Jacobi matrix.
An infinite matrix $J$ is called a Jacobi matrix if it has a tridiagonal form
\[
J =
\left(
\begin{array}{cccccc}
\beta_0 & \lambda_0 & 0 & 0 & 0 & \ldots \\
\lambda_0 & \beta_1 & \lambda_1 & 0 & 0 & \ldots\\
0 & \lambda_1 & \beta_2 & \lambda_2 &  0 & \ldots\\
0 & 0 &\lambda_2 & \beta_3 & \ddots & \\
0 & 0 & 0 & \ddots & \ddots & \\
\vdots & \vdots & \vdots & & & \\
\end{array}
\right),
\]
where the diagonal coefficients $\beta_n$  are real,
while the off-diagonal coefficients $\lambda_n$   are positive.
There exists a one to one correspondence between positive definite sequences $m_n$
and Jacobi matrices $J$ given by
\[
m_n \ = \ \iloczyn{J^n \delta_0}{\delta_0},
 \]
where $J$ is regarded as a symmetric unbounded operator on $\ell^2(\N)$.

\hspace{7pt}

Uniqueness of  the measure $\mu$ on the line turns out to be equivalent to
essential selfadjointness of $J$ on the subspace of finitely supported sequences in $\ell^2(\N)$.
Moreover, in the case when the moment problem $m_n$ is indeterminate,
a description of all the solutions $\mu$ is
related to a description of all selfadjoint extensions of $J.$

\hspace{7pt}

Selfadjointness of an unbounded  operator is an  important notion on more general grounds.
If a symmetric operator admits a selfadjoint extension, or even better, is essentially selfadjoint,
then the whole machinery of the spectral theory becomes available.

\hspace{7pt}

We take up the problem of essential selfadjointness of a Jacobi matrix on spaces
which are natural generalizations of $\ell^2(\N).$
The linear infinite tree $\N$ of nonnegative integers has two obvious extensions.
We may consider a homogeneous tree branching out from each vertex into a fixed number of edges
directed either downwards (the case of the tree $\Gamma$ with one origin)
or upwards (the case of tree $\Lambda$ with one end at infinity).
We consider a Jacobi matrix $J$ as a symmetric operator acting in
the space of all square-summable functions defined on the partially ordered set of vertices of these trees.
The domain of $J$ consists of finitely supported functions.
The main goal of this work is to investigate an  essential selfadjointness of the operator $J$.
In order to do this we look at its deficiency space.
It is described by  a recurrence relation, whose solutions yield systems of orthogonal polynomials.
Here the tools of the theory of orthogonal polynomials prove to be very useful.
It turns out that essential selfadjointness of $J$ depends on the structure of the tree
and, surprisingly, the behavior of $J$ is completely different in both cases.

\hspace{7pt}

In fact, the main result of Chapter 2 of the   paper states that

\hspace{7pt}

\noindent
\emph{
The matrix $J$ in the case of the tree $\Lambda$ is always essentially selfadjoint regardless of its coefficients.
Furthermore, $J$ has a pure point spectrum, i.e. there
is an orthonormal basis    consisting of eigenvectors for $J.$
}

\hspace{7pt}

In the case of $\Gamma,$ essential selfadjointness of $J$ depends on its projection on the one-dimensional tree $\N$.
Namely, we associate the Jacobi operator $J$ acting on the tree $\Gamma$
with some classical Jacobi matrix $J^r$ acting in $\ell^2(\N)$
which corresponds to the restriction of $J$ to the functions constant on levels of $\Gamma$.

\hspace{7pt}

The main result of Chapter 1 is

\hspace{7pt}

\noindent
\emph{
The operator $J$ is essentially selfadjoint if and only if $J^r$ is essentially selfadjoint.
}

\hspace{7pt}

One should not be misled by the apparent similarity to the classical  case.
The picture becomes clearer when we consider  the case when   $J$
is not essentially selfadjoint.
Then its deficiency space is much bigger than in the case of $\ell^2(\N)$
when it is just one dimensional.
We give a description of the nontrivial deficiency space of $J$ on the tree $\Gamma$.
It resembles the theory of harmonic functions since a Poisson-like kernel shows up there.
We prove that functions in the deficiency space are determined by their boundary values via
the Poisson integral. The spectral decomposition of selfadjoint extensions of $J$ is given explicitly.
In particular, we show that any such extension has a pure point spectrum.

\hspace{17pt}

\begin{center}
{\textsc{ACKNOWLEDGEMENT}}
\end{center}

\hspace{7pt}

The authors are grateful to Christian Berg and Pawe{\l} G\l owacki for
many helpful remarks which improved the exposition.


\clearpage

\begin{center}
{\textsc{PRELIMINARIES}}
\end{center}

\hspace{7pt}

\begin{center}
{\textsc{Selfadjoint extensions of  symmetric operators}}
\end{center}

\hspace{7pt}

Let $\h$ be a Hilbert space with  inner product $(\cdot,\cdot)$.
Let $A$ be a linear operator with domain $D(A) \subset \h$ which is dense in $\h$.

\bigskip

An operator $A$ is said to be \textbf{symmetric} if
\[
(Ax,y)=(x,Ay), \qquad x,y\in D(A).
\]

\bigskip

For an operator $A$ the \textbf{adjoint}  $A^*$ of $A$ is the operator defined on the space
\[
D(A^*)=\big\{x\in\h\colon\qquad\exists z\in\h\quad \forall y\in D(A)\quad (Ay,x)=(y,z) \big\}
\]
as follows: for $x\in D(A^*)$ we set $A^*x=z$.

\bigskip

If $A$ is symmetric and $D(A)=D(A^*)$, i.e. $A=A^*$, then $A$ is said to be \textbf{selfadjoint}.

\bigskip

If the graph of an operator $A$
\[
\big\{ \langle x,Ax\rangle\in\h\times \h\colon \quad x\in D(A) \big\}
\]
is a closed set then $A$ is called \textbf{closed}.

\bigskip

Clearly, each symmetric operator can be extended to a closed operator by taking a closure of its graph.
If this extension is already  selfadjoint then $A$ is called \textbf{essentially selfadjoint}.

\bigskip

For a  symmetric operator $A$ and a fixed complex number $z\notin\R$ we define the \textbf{deficiency space} of A by
\[
N_z\ =\ \Big(\obraz (A-\bar{z}I)\Big)^{\perp},
\]
where $^\perp$ denotes the orthogonal complement in $\h$.
It is known that the dimension of $N_z$ is constant on each of the half-planes $\im z >0$ and $\im z< 0$. These two numbers
$\dim N_i$ and $\dim N_{-i}$ are called the \textbf{deficiency indices} of $A$.

\bigskip

\begin{tw}
The deficiency space $N_z$ is a linear eigenspace of the operator $A^*$ associated with the eigenvalue $z$.
\end{tw}

\bigskip

\begin{tw}\label{admits}
A symmetric operator admits a selfadjoint extension if and only if its deficiency indices are equal.
\end{tw}

\bigskip

\begin{tw}\label{dodanie_ograniczonego}
Let $A$  be a symmetric operator and $B$ be a bounded selfadjoint operator.
Then the operators $A$ and $A+B$  have the same  deficiency indices.
\end{tw}

\bigskip

\begin{tw}\label{samosp_a_defekt}
A symmetric operator is essentially selfadjoint if and only if its deficiency space is trivial for any $z\notin \R$.
(i.e. its deficiency indices are zeros.)
\end{tw}

\bigskip

Presented facts can be found in many books, for instance in \cite{reedsimon1}, \cite{reedsimon2}, \cite{szwarc}.

\hspace{7pt}

\begin{center}
{\textsc{Classical Jacobi Matrices}}
\end{center}

\hspace{7pt}

A Jacobi matrix $J$, i.e. the matrix of the following form
\begin{equation}\label{matrix}
J =
\left(
\begin{array}{cccccc}
\beta_0 & \lambda_0 & 0 & 0 & 0 & \ldots \\
\lambda_0 & \beta_1 & \lambda_1 & 0 & 0 & \ldots\\
0 & \lambda_1 & \beta_2 & \lambda_2 &  0 & \ldots\\
0 & 0 &\lambda_2 & \beta_3 & \ddots & \\
0 & 0 & 0 & \ddots & \ddots & \\
\vdots & \vdots & \vdots & & & \\
\end{array}
\right),
\end{equation}
where $\beta_n$ are real and $\la_n$ are positive,
can be regarded as a linear operator in the Hilbert space $\ell^2(\N)$
with domain
$D(J)=\textrm{lin} \{\delta_0,\delta_1,\delta_2, \ldots\}$.
Then the action of $J$ on the characteristic function $\de_n$ of the point $n$ is expressed by the formula
\begin{equation}\label{rekurencja}
J \de_n \ = \ \la_{n-1}\de_{n-1} \ + \ \beta_n\delta_n \ + \ \la_n \delta_{n+1}, \qquad n\ge 0
\end{equation}
(we adopt the convention that $\la_{-1}=\de_{-1}=0$).

\bigskip

There are two sequences $p_n(x)$ and $q_n(x)$ of the orthogonal polynomials
associated with a Jacobi matrix $J$.
They are the solutions to the recurrence relation
\begin{equation}\label{a_n}
x\cdot a_n \ = \ \la_{n-1}a_{n-1} \ + \ \beta_n a_n \ + \ \la_n a_{n+1}, \qquad n\ge 1,
\end{equation}
with given initial conditions $a_0$ and $a_1$.
Taking  $a_0=1$ and $a_1=\frac{1}{\la_0}(x-\beta_0)$ gives $a_n=p_n$;
while $a_0=0$ and $a_1=\frac{1}{\la_0}$ give $a_n=q_n$.
It is known that all roots of these polynomials are real (see e.g. \cite{chihara}).

\bigskip

The following facts concerning basic properties of Jacobi matrices can be found in many books and articles,
for instance, in \cite{akhiezer}, \cite{berg}, \cite{chihara}, \cite{simon}, \cite{szwarc}.

\bigskip

From (\ref{matrix}) we can see that the operator $J$ is symmetric.
In view of Theorem \ref{admits} the following one implies that $J$ has a selfadjoint extension.

\bigskip

\begin{tw}\label{jednowymiarowy_defekt}
The deficiency indices of the operator $J$ are equal either to (0,0) or to (1,1).
In the former case $J$ is essentially selfadjoint.
In the latter case a selfadjoint extension of $J$ is not unique.
\end{tw}

\bigskip

\begin{tw}{\bf{(The Hamburger criterion)}}  \label{hamburger}
A Jacobi matrix $J$ is essentially selfadjoint if and only if at least one of
the series $\sum p_n(0)^2$ and $\sum q_n(0)^2$ is divergent.
\end{tw}

\bigskip

\begin{tw}\label{n_ekstremalne}
Let $\tilde{J}$ be a selfadjoint extension of $J$ in the indeterminate case and $E(x)$ be the resolution of the identity associated
with~$\tilde{J}$.
Then the support of the measure $$d\sigma(x)=d(E(x)\de_0,\de_0)$$ is a discrete set
and coincides with the spectrum of the operator $\tilde{J}$.
\end{tw}

\bigskip

The concept of selfadjointness of the operator $J$ is important in the theory of the classical orthogonal polynomials.
The measure which is the solution to the moment problem $m_n=(J^n\delta_0,\delta_0)$
is unique if and only if the Jacobi matrix $J$ is essentially selfadjoint.

\hspace{7pt}

\begin{center}
{\textsc{Jacobi Matrices on homogeneous trees}}
\end{center}

\hspace{7pt}

The set of non-negative integers  $\N$ can be identified with a linear infinite tree
with a \emph{natural} order.
\bigskip
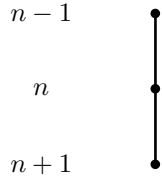
\footnotesize
\begin{center}
\begin{pspicture}(1,0)(3,2)
\psline[linewidth=1pt](3,2)(3,0)
\psdots[dotstyle=*](3,0)(3,2)(3,1)
\rput(1.5,0){$n+1$}
\rput(1.5,1){$n$}
\rput(1.5,2){$n-1$}
\end{pspicture}
\end{center}
\normalsize
\bigskip

It seems that there are two natural generalizations of this configuration: from each vertex
there is a fixed number (greater than 1) of edges either pointing downward
(a tree with one origin) or upward (a tree with one end).
\bigskip
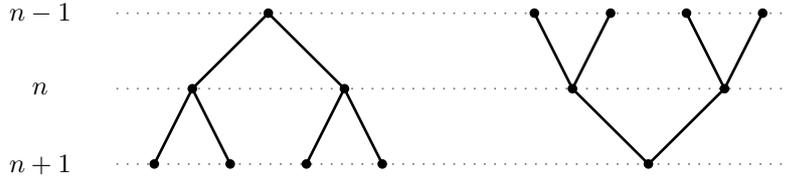
\footnotesize
\begin{center}
\begin{pspicture}(1,0)(12,2)
\psline[linestyle=dotted,linecolor=gray](2.5,2)(11.5,2)
\psline[linestyle=dotted,linecolor=gray](2.5,1)(11.5,1)
\psline[linestyle=dotted,linecolor=gray](2.5,0)(11.5,0)
\psline[linewidth=1pt](3,0)(3.5,1)(4.5,2)(5.5,1)(6,0)
\psline[linewidth=1pt](3.5,1)(4,0)
\psline[linewidth=1pt](5,0)(5.5,1)
\psline[linewidth=1pt](8,2)(8.5,1)(9.5,0)(10.5,1)(11,2)
\psline[linewidth=1pt](8.5,1)(9,2)
\psline[linewidth=1pt](10,2)(10.5,1)
\psdots[dotstyle=*](3,0)(4,0)(5,0)(6,0)(3.5,1)(5.5,1)(4.5,2)
\psdots[dotstyle=*](9.5,0)(8.5,1)(10.5,1)(8,2)(9,2)(10,2)(11,2)
\rput(1.5,0){$n+1$}
\rput(1.5,1){$n$}
\rput(1.5,2){$n-1$}
\end{pspicture}
\end{center}
\normalsize

\hspace{7pt}


\newpage

\section{\textsc{A JACOBI OPERATOR ON A TREE WITH ONE ORIGIN}}

\hspace{7pt}

For a fixed number $d\in\{2,3,4,\ldots\}$ we consider an infinite homogeneous tree of degree  $d$,
i.e. an infinite connected graph with a distinguished vertex (root) $e$ and a partial order
such that  for each vertex there are exactly $d$ edges downward (with respect to $e$).
In other words, each vertex has $d$ successors and each vertex different from the root
has one predecessor.
For a fixed vertex $x$ let $x_0$ and $x_i$ ($i=1,2,\ldots,d$) denote respectively the predecessor and the successors of $x$.

\bigskip

For instance, if $d=3$, the top levels of the tree look as follows

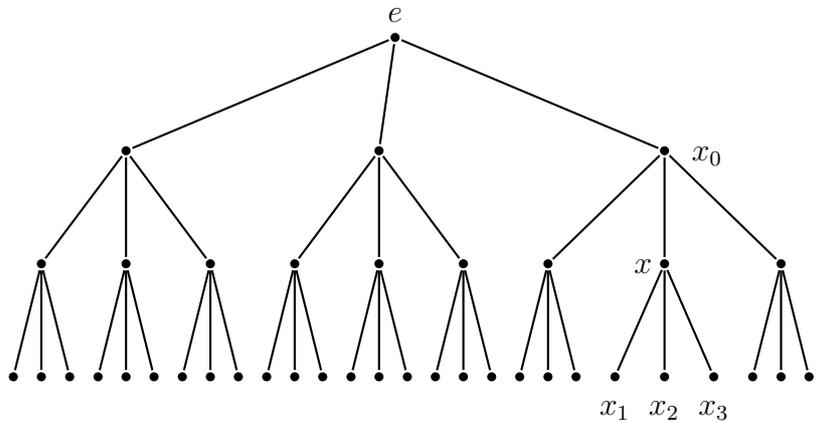
\vspace{10pt}

\begin{center}
\pstree[nodesep=0.3mm,treemode=D,treefit=loose,treesep=0.25,levelsep=1.5cm]
{\Tc*{0.6mm}~[tnpos=a,tnsep=0mm]{$e$}}
{

\pstree
{\Tc*{0.6mm}~[tnpos=l,tnsep=3mm]{$  $}}
{
  \pstree
  {\Tc*{0.6mm}~[tnpos=l,tnsep=1mm]{$ $}}
  {
   \Tc*{0.6mm}~[tnpos=b,tnsep=1mm]{$ $}
   \Tc*{0.6mm}~[tnpos=b,tnsep=1mm]{$ $}
   \Tc*{0.6mm}~[tnpos=b,tnsep=1mm]{$ $}
  }
    \pstree
     {\Tc*{0.6mm}~[tnpos=l,tnsep=1mm]{$ $}}
     {  \Tc*{0.6mm}~[tnpos=b,tnsep=1mm]{$ $}
        \Tc*{0.6mm}~[tnpos=b,tnsep=1mm]{$ $}
        \Tc*{0.6mm}~[tnpos=b,tnsep=1mm]{$ $}
     }
       \pstree
        {\Tc*{0.6mm}~[tnpos=l,tnsep=1mm]{$ $}}
        {
        \Tc*{0.6mm}~[tnpos=b,tnsep=1mm]{$ $}
        \Tc*{0.6mm}~[tnpos=b,tnsep=1mm]{$ $}
        \Tc*{0.6mm}~[tnpos=b,tnsep=1mm]{$ $}
  }
}

 \pstree
{\Tc*{0.6mm}~[tnpos=l,tnsep=3mm]{$  $}}
{
  \pstree
  {\Tc*{0.6mm}~[tnpos=l,tnsep=1mm]{$ $}}
  {
   \Tc*{0.6mm}~[tnpos=b,tnsep=1mm]{$ $}
   \Tc*{0.6mm}~[tnpos=b,tnsep=1mm]{$ $}
   \Tc*{0.6mm}~[tnpos=b,tnsep=1mm]{$ $}
  }
    \pstree
     {\Tc*{0.6mm}~[tnpos=l,tnsep=1mm]{$ $}}
     {  \Tc*{0.6mm}~[tnpos=b,tnsep=1mm]{$ $}
        \Tc*{0.6mm}~[tnpos=b,tnsep=1mm]{$ $}
        \Tc*{0.6mm}~[tnpos=b,tnsep=1mm]{$ $}
     }
       \pstree
        {\Tc*{0.6mm}~[tnpos=l,tnsep=1mm]{$ $}}
        {
        \Tc*{0.6mm}~[tnpos=b,tnsep=1mm]{$ $}
        \Tc*{0.6mm}~[tnpos=b,tnsep=1mm]{$ $}
        \Tc*{0.6mm}~[tnpos=b,tnsep=1mm]{$ $}
  }
}

\pstree
{\Tc*{0.6mm}~[tnpos=r,tnsep=3mm]{${x}_0$}}
{
  \pstree
  {\Tc*{0.6mm}~[tnpos=l,tnsep=1mm]{$ $}}
  {
   \Tc*{0.6mm}~[tnpos=b,tnsep=1mm]{$ $}
   \Tc*{0.6mm}~[tnpos=b,tnsep=1mm]{$ $}
   \Tc*{0.6mm}~[tnpos=b,tnsep=1mm]{$ $}
  }
    \pstree
     {\Tc*{0.6mm}~[tnpos=l,tnsep=1mm]{$x$}}
     {  \Tc*{0.6mm}~[tnpos=b,tnsep=1mm]{${x}_1$}
        \Tc*{0.6mm}~[tnpos=b,tnsep=1mm]{${x}_2$}
        \Tc*{0.6mm}~[tnpos=b,tnsep=1mm]{${x}_3$}
     }
       \pstree
        {\Tc*{0.6mm}~[tnpos=l,tnsep=1mm]{$ $}}
        {
        \Tc*{0.6mm}~[tnpos=b,tnsep=1mm]{$ $}
        \Tc*{0.6mm}~[tnpos=b,tnsep=1mm]{$ $}
        \Tc*{0.6mm}~[tnpos=r,tnsep=1mm]{$.$}
  }
}

}
\end{center}

\vspace{10pt}

The set of all vertices of the tree will be denoted by $\Gamma_d$.
There is a  natural distance $\dist(\cdot, \cdot)$ in  $\gd$ counting the number of edges in the unique path connecting
two fixed vertices.
The {\emph{length}} of a vertex is, by definition, its distance from the root $e$, i.e.
\[
|x| \ =\ \dist(x,e).
\]
For the vertex $x$ marked in the picture above we have $|x|=2$ and $|x_1|=|{x}_2|=3$.

\bigskip

The space $\ell^2(\gd)$ of all square-summable functions on $\gd$, i.e.
\[
\ell^2(\gd)\ = \
\big\{f\in\C^{\gd}\ \colon\quad\sum\limits_{x\in\gd}|f(x)|^2<\infty \big\},
\]
is a Hilbert space with the standard inner product
\[
(f,g)\ =\ \sum\limits_{x\in\gd} f(x)\overline{g(x)}.
\]
We write $\delta_x$ for the characteristic function of  the one point set $\{x\}$.
Let $\F$ denote the space of all functions with finite support:
\[
\F\ = \ \lin \{\ \delta_x\colon\ x\in\gd \ \}.
\]

\bigskip

Let $\lambda_0, \lambda_1, \lambda_2, \ldots$ be fixed positive numbers
and $\beta_0,\beta_1,\beta_2,\ldots$ be fixed real numbers.
We consider the Jacobi operator $J$ with domain
$$
D(J)\ =\ \F
\ \subset \ \ell^2(\gd),
$$
which acts as follows
\begin{equation}\label{definicja}
\begin{split}
J \de_e \ & = \qquad\qquad\qquad\ \beta_0\cdot\de_e \ + \ \la_0\cdot \big(\de_{e_1}+\ldots+\de_{e_d}\big)\\
J \delta_x \ & = \ \lambda_{n-1}\cdot\delta_{{x}_0}\ +\ \beta_n\cdot\delta_x
             \  +\ \lambda_n \cdot \big(\delta_{{x}_1} + \ldots +\delta_{{x}_d}\big), \quad n\ge 1,
\end{split}
\end{equation}
where $n=|x|$.
We adopt the convention that $\lambda_{-1}=\delta_{{e}_0}=0$. Then the action of $J$ can be expressed by the latter formula for all $n\ge0$.

\bigskip

It is elementary  that thus defined $J$  is a symmetric operator.

\bigskip
\bigskip

\begin{f}
The deficiency space $N_z(J)$ of the operator $J$ on $\ell^2(\gd)$ consists of all square-summable functions $v(x)$ on $\gd$
satisfying
\begin{equation}\label{rownania_defektu}
z v(x)=\lambda_{n-1}v({x}_0)+ \beta_nv(x) + \lambda_n \big(v({{x}_1}) + \ldots +v({{x}_d})\big)
\end{equation}
for all $|x|=n$ and  all $n\geq 0$.
\end{f}

\bigskip

\begin{dowod}
A function $v\in \ell^2(\gd)$ is orthogonal to the image $\obraz(J-\bar{z}I)$ if and only if
for each vertex $|x|=n$
\begin{eqnarray*}
0 &=&(v\;,\; (J-\bar{z}) \de_x)  \\
  &=&(v\;,\; \lambda_{n-1}\de_{{x}_o} + \beta_n\de_x
               + \lambda_n \big(\de_{{x}_1} + \ldots +\de_{{x}_d}\big)
               -\bar{z}\de_x)  \\
  &=& \lambda_{n-1}v({x}_0)+ \beta_nv(x)
               + \lambda_n \big(v({{x}_1}) + \ldots +v({{x}_d})\big)
               -z v(x).
\end{eqnarray*}
\end{dowod}

\bigskip
\bigskip

\noindent
{\textbf{Remark.}}
Although the domain of $J$ consists of functions with  finite support, note that the formula for
$J$  can be actually applied to  any  function on $\gd$.
Therefore we can write
\[
N_z(J)\ =\ \Big\{\, v\in\ldwa\colon \ \ Jv\, (x)\ =\ z \cdot v(x),\ \ \ x\in\gd \,\Big\}.
\]


\hspace{7pt}

\begin{center}
{\textsc{The one dimensional operator}}
\end{center}

\hspace{7pt}

We call a function on $\Gamma_d$ \emph{radial} if it is constant on each {\emph{level}} of $\gd$,
that is to say,  on each set of vertices of fixed length.
We will denote by $\ell^2_r (\Gamma_d) $
the space of all square-summable functions on $\gd$ which are radial.
Let $\chi_n$ denote the characteristic function of the  $n$th level.
Note that the normalized functions
\[
\mu_n (x)\ = \ \big(\sqrt{d}\,\big)^{-n} \cdot \chi_n(x) \ = \
\left\{ \begin{array}{ll}
d^{-n\slash 2} & \textrm{for }\ |x|=n,\\
0 & \textrm{for }\ |x| \neq n ,
\end{array} \right.
\]
form an orthonormal basis of $\ell^2_r (\gd)$.

\bigskip

Obviously,
\[
\chi_n \ = \ \sum\limits_{|x|=n} \de_x.
\]
Each vertex of length $n-1$ is a predecessor of exactly $d$ vertices of length $n$.
Therefore applying $J$ to the characteristic function $\chi_n$ gives
\[
J\chi_n \ =\ d\cdot\la_{n-1}\cdot\chi_{n-1} + \beta_n\cdot\chi_n
             +\la_n\cdot\chi_{n+1}.
\]
Since $\chi_n=\big(\sqrt{d}\,\big)^{n}\mu_n$, we have
\[
\big(\sqrt{d}\,\big)^{n} J\mu_n \ =\ d\big(\sqrt{d}\,\big)^{n-1}\la_{n-1}\mu_{n-1} + \big(\sqrt{d}\,\big)^{n}\beta_n\mu_n
             +\big(\sqrt{d}\,\big)^{n+1}\la_n\mu_{n+1}.
\]
The restriction of $J$ to the subspace $\ell^2_r(\gd)$ of radial functions will be denoted by $J^r$.
We thus  have
\[
D(J^r)\ = \ \lin \{\mu_0,\mu_1, \mu_2, \ldots \}\ \subset\ \ell^2_r(\gd)
\]
and
\begin{equation}\label{recurrence_j_r}
J^r\mu_n \ =\ \sqrt{d}\la_{n-1}\cdot\mu_{n-1} + \beta_n\cdot\mu_n
               + \sqrt{d}\la_n\cdot\mu_{n+1},\quad n\ge 0.
\end{equation}
In other words, we can identify the action of the radial operator $J^r$ with the matrix
\begin{equation}\label{matrix_j_r}
J^r \ = \
\left(
\begin{array}{cccccc}
\beta_0   & \pd\la_0 & 0 & 0 & 0  & \ldots\\
\pd\la_0  & \beta_1   & \pd\la_1 & 0 & 0 & \ldots\\
0 & \pd\la_1  & \beta_2   & \pd\la_2 &  0 & \ldots\\
0 & 0 & \pd\la_2  & \beta_3   & \pd\la_3 &   \ldots\\
0 & 0 & 0 & \pd\la_3 & \beta_4 & \ddots \\
\vdots & \vdots & \vdots & \vdots & \ddots & \ddots
\end{array}
\right).
\end{equation}
It means that $J^r$ on $\ell^2_r(\gd)$ can be regarded as a classical one dimensional Jacobi operator on $\ell^2(\N)$.
In particular, by Theorem \ref{jednowymiarowy_defekt}, its deficiency space $N_z(J^r)$ is either one dimensional or trivial.

\bigskip
\bigskip

\begin{f}\label{def}
A function $v\in \ell^2_r(\gd)$ belongs to the deficiency space  $N_z(J^r)$ of the radial operator $J^r$ on $\,\ell^2_r(\gd)$
if and only if
\begin{equation}\label{defekt_r}
z v(x)=\lambda_{n-1}v({x}_0)+ \beta_nv(x) + \lambda_n \big(v({{x}_1}) + \ldots +v({{x}_d})\big)
\end{equation}
for each $n\ge 0$ and for each $|x|=n$. Moreover,
\[
N_z(J^r)\ \subseteq \ N_z(J).
\]
\end{f}

\bigskip

\begin{dowod}
Let a function $v\in \ell^2_r(\gd)$ be orthogonal to $\obraz(J^r-\bar{z}I)$, i.e.
\[
0 \ =\ \big(\,v\;,\; (J^r-\bar{z}) \chi_n\,\big), \qquad n\ge 0.
\]
We calculate
\begin{align*}
&(v,\ \ (J^r-\bar{z}) \chi_n) \\
 &=(v\;,\; d\cdot\lambda_{n-1}\cdot\chi_{n-1} + \beta_n\cdot\chi_n
               + \lambda_n\cdot \chi_{n+1} - \bar{z}\chi_n)  \\
 &=d\lambda_{n-1}\sum\limits_{|x|=n-1}v(x)+ \beta_n\sum\limits_{|x|=n}v(x)
               + \lambda_n \sum\limits_{|x|=n+1}v(x) - z \sum\limits_{|x|=n}v(x).
\end{align*}
Since $v$ is radial, we obtain
\[
0
\;=\;
d\lambda_{n-1}\cdot d^{n-1}v({x}_0)+ \beta_n\cdot d^n v(x)
               + \lambda_n\cdot d^{n+1} v({{x}_1}) - z\cdot d^n v(x).
\]
It follows that
\[
0
\;=\;
\lambda_{n-1}v({x}_0)+ \beta_n v(x)
               + \lambda_n d\cdot  v({{x}_1}) - z v(x)
\]
for each vertex $|x|=n$.
\end{dowod}

\bigskip
\bigskip

It turns out that the problem of the essential selfadjointness of $J$  on $\ell^2(\gd)$ can be reduced to the same problem
for the corresponding one dimensional operator $J^r$  on $\ell^2(\N)$.

\bigskip

\begin{tw}\label{iff}
The operator $J$  on $\ell^2(\gd)$ is essentially selfadjoint if and only if the corresponding one dimensional operator $J^r$ on $\ell^2(\N)$
is essentially selfadjoint.
\end{tw}

\bigskip

\begin{dowod}
By Theorem \ref{samosp_a_defekt} and Fact \ref{def}, it suffices to show that if $J$ is not essentially
selfadjoint, neither is the matrix $J^r$.
To this end, assume that the deficiency space $N_z(J)$ is nontrivial, i.e. there exists
\[
0 \ \neq \ f\ \in\  N_z(J).
\]
We will construct  a special function in a deficiency space.
This will allow us to show that $J^r$ is not essentially selfadjoint.

\bigskip

Let $x$ be a vertex in the support of $f$ of minimal length, i.e.
\[
f(x)\neq 0 \qquad\text{ and }\qquad f(y)=0\  \text{ for }\ |y|<|x|.
\]
Let  $\Gamma_x$ denote the subtree of the tree $\gd$ which has its root at $x$.
The subtree $\Gamma_x$ is marked in the picture below.

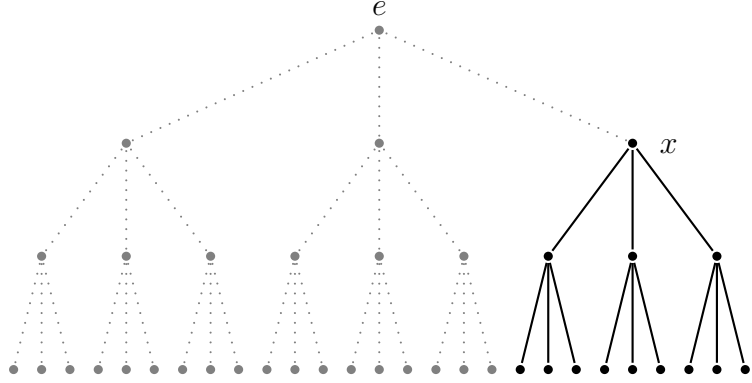
\begin{center}
\pstree[nodesep=0.3mm,treemode=D,treefit=loose,treesep=0.25,levelsep=1.5cm,linestyle=dotted,linecolor=gray]
{\Tc*[linecolor=gray]{0.6mm}~[tnpos=a,tnsep=0mm]{$e$}}
{

\pstree
{\Tc*{0.6mm}~[tnpos=l,tnsep=3mm]{$  $}}
{
  \pstree
  {\Tc*{0.6mm}~[tnpos=l,tnsep=1mm]{$ $}}
  {
   \Tc*{0.6mm}~[tnpos=b,tnsep=1mm]{$ $}
   \Tc*{0.6mm}~[tnpos=b,tnsep=1mm]{$ $}
   \Tc*{0.6mm}~[tnpos=b,tnsep=1mm]{$ $}
  }
    \pstree
     {\Tc*{0.6mm}~[tnpos=l,tnsep=1mm]{$ $}}
     {  \Tc*{0.6mm}~[tnpos=b,tnsep=1mm]{$ $}
        \Tc*{0.6mm}~[tnpos=b,tnsep=1mm]{$ $}
        \Tc*{0.6mm}~[tnpos=b,tnsep=1mm]{$ $}
     }
       \pstree
        {\Tc*{0.6mm}~[tnpos=l,tnsep=1mm]{$ $}}
        {
        \Tc*{0.6mm}~[tnpos=b,tnsep=1mm]{$ $}
        \Tc*{0.6mm}~[tnpos=b,tnsep=1mm]{$ $}
        \Tc*{0.6mm}~[tnpos=b,tnsep=1mm]{$ $}
  }
}

 \pstree
{\Tc*{0.6mm}~[tnpos=l,tnsep=3mm]{$  $}}
{
  \pstree
  {\Tc*{0.6mm}~[tnpos=l,tnsep=1mm]{$ $}}
  {
   \Tc*{0.6mm}~[tnpos=b,tnsep=1mm]{$ $}
   \Tc*{0.6mm}~[tnpos=b,tnsep=1mm]{$ $}
   \Tc*{0.6mm}~[tnpos=b,tnsep=1mm]{$ $}
  }
    \pstree
     {\Tc*{0.6mm}~[tnpos=l,tnsep=1mm]{$ $}}
     {  \Tc*{0.6mm}~[tnpos=b,tnsep=1mm]{$ $}
        \Tc*{0.6mm}~[tnpos=b,tnsep=1mm]{$ $}
        \Tc*{0.6mm}~[tnpos=b,tnsep=1mm]{$ $}
     }
       \pstree
        {\Tc*{0.6mm}~[tnpos=l,tnsep=1mm]{$ $}}
        {
        \Tc*{0.6mm}~[tnpos=b,tnsep=1mm]{$ $}
        \Tc*{0.6mm}~[tnpos=b,tnsep=1mm]{$ $}
        \Tc*{0.6mm}~[tnpos=b,tnsep=1mm]{$ $}
  }
}

\pstree[linecolor=black,linestyle=solid]
{\Tc*[linecolor=black]{0.6mm}~[tnpos=r,tnsep=3mm]{$ x $}}
{
  \pstree
  {\Tc*{0.6mm}~[tnpos=l,tnsep=1mm]{$ $}}
  {
   \Tc*{0.6mm}~[tnpos=b,tnsep=1mm]{$ $}
   \Tc*{0.6mm}~[tnpos=b,tnsep=1mm]{$ $}
   \Tc*{0.6mm}~[tnpos=b,tnsep=1mm]{$ $}
  }
    \pstree
     {\Tc*{0.6mm}~[tnpos=l,tnsep=1mm]{$ $}}
     {  \Tc*{0.6mm}~[tnpos=b,tnsep=1mm]{$ $}
        \Tc*{0.6mm}~[tnpos=b,tnsep=1mm]{$ $}
        \Tc*{0.6mm}~[tnpos=b,tnsep=1mm]{$ $}
     }
       \pstree
        {\Tc*{0.6mm}~[tnpos=l,tnsep=1mm]{$ $}}
        {
        \Tc*{0.6mm}~[tnpos=b,tnsep=1mm]{$ $}
        \Tc*{0.6mm}~[tnpos=b,tnsep=1mm]{$ $}
        \Tc*{0.6mm}~[tnpos=b,tnsep=1mm]{$ $}
  }
}

}
\end{center}

\vspace{10pt}

In the proof we are going to apply an averaging operator $E$.

\bigskip

\begin{lem}\label{kontrakcja}
In $\ell^2(\gd)$ the averaging operator
\begin{equation}\label{E}
E f (w) \ = \ \frac{1}{d^{{\scriptscriptstyle |w|}}} \ \sum\limits_{|y|=|w|} f(y)
\end{equation}
is a selfadjoint projection.
\end{lem}

\bigskip

\begin{dowod}
For any functions $f,g\in \F$ we have
\begin{align*}
(\; Ef\;,\;g\;)\ & =\ \sum\limits_{w\in\gd}Ef(w)\cdot \overline{g(w)}
\ =\ \sum\limits_{k=0}^{\infty} \sum\limits_{|w|=k}\Bigg( \frac{1}{d^k}  \sum\limits_{|y|=k} f(y) \Bigg)\cdot \overline{g(w)} \\
& =\ \sum\limits_{k=0}^{\infty} \frac{1}{d^k}\sum\limits_{|w|=k} \sum\limits_{|y|=k} f(y)  \overline{g(w)}.
\end{align*}
Reversing the order of summation yields
\begin{align*}
(\; Ef\;,\;g\;)\ &  =\ \sum\limits_{k=0}^{\infty} \sum\limits_{|y|=k}\Bigg( \frac{1}{d^k}  \sum\limits_{|w|=k} \overline{g(w)} \Bigg)\cdot f(y) \\
& =\ \sum\limits_{w\in\gd}f(w)\cdot \overline{Eg(w)}
\ =\ (\; f\;,\;Eg\;),
\end{align*}
which proves the symmetry.

Now, by the Schwarz inequality,
\begin{eqnarray*}
\| E f \|^2
&=&   \sum\limits_{w\in\gd} \big|Ef(w)\big|^2
      \ =\   \sum\limits_{k=0}^{\infty}\ \sum\limits_{|w|=k}\ \Big|\, d^{-k} \sum\limits_{|y|=k} f(y) \,\Big|^{\,2}\\
&\le& \sum\limits_{k=0}^{\infty}\ d^k \ d^{-2k} \ \bigg(\sum\limits_{|y|=k} |f(y)| \bigg)^2\\
&\le& \sum\limits_{k=0}^{\infty}\ d^{-k} \ \bigg(\sum\limits_{|y|=k} 1^2 \bigg)\ \bigg(\sum\limits_{|y|=k} |f(y)|^2\bigg)\\
&=&   \sum\limits_{k=0}^{\infty}\ \sum\limits_{|y|=k}\ \big|f(y)\big|^2
\ = \ \sum\limits_{y\in\gd} \big|f(y)\big|^2 \ = \ \|f\|^2,
\end{eqnarray*}
whence the norm of $E$ is less than or  equal to 1.
Moreover, for a radial function $f\in\ell^2_r(\gd)$ we obtain the equality $\|Ef\|=\|f\|$.
\end{dowod}

\bigskip

We denote by $f_x$ the restriction of $f$ to the subtree $\Gamma_x$.
Let $k=|x|$.
The symbol $E_x$ will denote the averaging operator on $\Gamma_x$.
More precisely, $E_x(g)$ is the mean value
of a function $g$ on each level of the subtree $\Gamma_x$ around its root $x$:
\[
E_x\ \colon\quad \ell^2(\Gamma_x)\ \longrightarrow \ \ell^2_r(\Gamma_x)
\]
and
\begin{equation}\label{E_x}
E_x\, g\,(y)\ =\ d^{-(|y|-k)}\cdot \sum\limits_{\stackrel{t\in\Gamma_x}{|t|=|y|}} g(t).
\end{equation}
By Lemma \ref{kontrakcja}, it is obvious that $E_x$ is a contraction on $\ell^2(\Gamma_x)$.
Thus the function $E_x(f_x)$ is square-summable and radial on $\Gamma_x$.
Since neither restricting nor taking the mean value around $x$  change the value at $x$,
the function $E_x(f_x)$ is nonzero.
In order to belong to a deficiency space it needs to satisfy appropriate equations.
Since $f$, as an element of the deficiency space $N_z(J)$,  satisfies all the recurrence equations (\ref{rownania_defektu}),
its restriction $f_x$ satisfies those of them which are related to the restriction of $J$  to $\Gamma_x$.
Indeed, in each vertex of $\Gamma_x$ different from $x$  the equations and  values  remain unchanged.
Therefore,  only the equation in $x$ can raise doubts.
However, at $x$ we have
\[
z f_x(x)\;=\; 0 + \beta_k f_x(x)
               + \lambda_k \big(f_x({{x}_1}) + \ldots +f_x({{x}_d})\big),
\]
which is consistent with the convention in (\ref{definicja})
applied to the operator $J$ with coefficients shifted by $k$.
The corresponding radial operator  is expressed by the matrix
\begin{equation}\label{macierz_obcieta}
J^r_k \ =\
\left(
\begin{array}{cccccc}
\beta_k   & \pd\la_k & 0 & 0 & 0  & \ldots\\
\pd\la_k  & \beta_{k+1}   & \pd\la_{k+1} & 0 & 0 & \ldots\\
0 & \pd\la_{k+1}  & \beta_{k+2}   & \pd\la_{k+2} &  0 & \ldots\\
0 & 0 & \pd\la_{k+2}  & \beta_{k+3}   & \pd\la_{k+3} &   \ldots\\
0  &  0 & 0 & \pd \la_{k+3} & \beta_{k+4} & \ddots \\
\vdots & \vdots & \vdots & \vdots & \ddots & \ddots
\end{array}
\right).
\end{equation}
It is immediate that taking the mean value on levels does not
affect the recurrence relation described above.
Hence
\[
0\ \neq\  E_x(f_x)\  \in\  N_z(J^r_k),
\]
i.e. the matrix $J^r_k$ is not essentially selfadjoint.
We add to $J^r_k$ an extra first column and first row
consisting on zeros.
We also add an extra first coordinate with value zero to the vector $E_x(f_x)$.
We thus get one additional equation in the description of the deficiency space of the new operator (cf. (\ref{defekt_r}))
which is obviously satisfied  as the equation consists of zeros.
Hence, also the extended matrix is not essentially selfadjoint.
Therefore, the matrix with exactly $k$ extra zero columns and rows
\[
\left(
\begin{array}{ccccc}
0 & 0   & 0& 0&\cdots\\
0 & 0   & 0& 0&\cdots\\
0 & 0  &  & & \\
0 &0  &   &J_k^r &\\
\vdots  & \vdots &&&
\end{array}
\right)
\]
is not essentially selfadjoint. Next we add to it a symmetric  finite dimensional operator of the form
\[
\left(
\begin{array}{cccccc}
\beta_0   & \pd\la_0 & 0 & 0 & 0  & \ldots\\
\pd\la_0  & \ddots   & \ddots &  & \vdots & \\
0 & \ddots  & \beta_{k-1}   & \pd\la_{k-1} &  0 & \\
0 &  & \pd\la_{k-1}  & 0   & 0 &  \\
0 & \ldots & 0 & 0 & 0 &  \ddots\\
\vdots &  &  & & \ddots & \ddots
\end{array}
\right).
\]
Since it is selfadjoint and bounded, the operator $J^r$, by Theorem \ref{dodanie_ograniczonego}, is not essentially selfadjoint.
\end{dowod}

\bigskip
\bigskip

\noindent
\textbf{Remark}.
We have associated with $J$ in $\ldwa$ the radial operator $J^r$ acting in $\ell^2_r(\gd)$, which can be identified with $\ell^2(\N)$.
These two matrices
\[
J  =
\left(
\begin{array}{cccc}
\beta_0   & \la_0 & 0 &  0\\
\la_0  & \beta_1   & \la_1 &  0\\
0 & \la_1  & \beta_2   &  \ddots\\
0 & 0 & \ddots & \ddots
\end{array}
\right)\ \ \textrm{ and } \ \
J^r  =
\left(
\begin{array}{cccc}
\beta_0   & \pd\la_0 & 0 &  0\\
\pd\la_0  & \beta_1   & \pd\la_1 &  0\\
0 & \pd\la_1  & \beta_2   &  \ddots\\
0 & 0 & \ddots & \ddots
\end{array}
\right)
\]
do not have to be essentially selfadjoint at the same time.
Let us consider an example. For $d=2$ let $\beta_n=\la_n+\la_{n-1}$ and $\beta_0=\la_0$. Then
\[
J\ =\
\left(
\begin{array}{ccccc}
\la_0   & \la_0 & & &  \\
\la_0  & \la_0+\la_1   & \la_1 & &  \\
 & \la_1  & \la_1+\la_2 &  \la_2 &    \\
 & & \la_2 & \la_2+\la_3 &\ddots \\
 &&&\ddots & \ddots
\end{array}
\right).
\]
The recurrence relation associated with $J$ (cf. (\ref{a_n})) is of the following form
\begin{align*}
xa_n \  
&= \ \la_{n-1} a_{n-1}+(\la_n+\la_{n-1})a_n +\la_n a_{n+1}\\
& = \ (a_{n-1}+a_n)\la_{n-1}+(a_n+a_{n+1})\la_{n}, \qquad n\ge 1.
\end{align*}
In particular, for $x=0$ we get
\begin{align*}
a_{n+1}(0)\ &= \ -\,\frac{\la_{n-1}}{\la_n}\big(a_{n-1}(0)+a_n(0)\big) - a_n(0).
\end{align*}
For the sequence $p_n(0)$ (cf. (\ref{a_n})) we get $a_0(0)=p_0(0)=1$ and $a_1(0)=p_1(0)=-1$.
Consequently, by induction $p_n(0)=(-1)^n$.
Hence the series $\sum p_n(0)^2$ is divergent. By the Hamburger criterion (cf. Theorem \ref{hamburger}),
the matrix $J$ is essentially selfadjoint.

The corresponding matrix on the tree $\Gamma_2$  is of the form (cf. \ref{matrix_j_r})
\[
J^r\ =\
\left(
\begin{array}{ccccc}
\la_0   & \sqrt{2}\la_0 & & & \\
\sqrt{2}\la_0  & \la_0+\la_1   & \sqrt{2}\la_1 & &  \\
 & \sqrt{2}\la_1  & \la_1+\la_2  & \sqrt{2}\la_2  &    \\
 & & \sqrt{2}\la_2 &  \la_2+\la_3  & \ddots \\
 &&&\ddots & \ddots
\end{array}
\right).
\]
Let $\la_n=2^n$. Then $\la_{n-1}+\la_n =  3\cdot 2^{n-1}$. Hence for $x=0$ the general solution  to the recurrence relation
\[
\sqrt{2}\cdot a_{n-1} +3 \cdot a_n +2\sqrt{2}\cdot a_{n+1} =0, \qquad n\ge 1,
\]
is
\[
a_n\ = \ \Big(\frac{1}{\sqrt{2}}\Big)^n \big(\,c_1 \cdot \cos n\theta + c_2 \cdot \sin n\theta \,\big).
\]
Thus the series
\[
\sum |a_n|^2 \ \le \ \big( |c_1|^2 +|c_2|^2 \big)\, \sum \frac{1}{2^n}
\]
is always convergent. It means that both series  $\sum p_n(0)^2$ and $\sum q_n(0)^2$
(cf. (\ref{a_n}))
are convergent.  By the Hamburger criterion (cf. Theorem \ref{hamburger}), the matrix  $J^r$ is not essentially selfadjoint.
\bigskip


\hspace{7pt}

\begin{center}
{\textsc{The description of the deficiency space}}
\end{center}

\hspace{7pt}

We are going to write down the nontrivial deficiency space $N_z(J)$ as a sum
of spaces associated with vertices of $\gd$.

\bigskip

Fix a vertex $x$ of length $k$. Let $J_k$ denote the truncated matrix
\[
J_k\ = \
\left(
\begin{array}{cccccc}
\beta_k   & \la_k & 0 & 0 & 0  & \ldots\\
\la_k  & \beta_{k+1}   & \la_{k+1} & 0 & 0 & \ldots\\
0 & \la_{k+1}  & \beta_{k+2}   & \la_{k+2} &  0 & \ldots\\
0 & 0 & \la_{k+2}  & \beta_{k+3}   & \la_{k+3} &   \ldots\\
0 & 0 & 0 & \la_{k+3} & \beta_{k+4} & \ddots \\
\vdots & \vdots & \vdots & \vdots & \ddots & \ddots
\end{array}
\right).
\]

Observe that $\Gamma_x$ which is the subtree of $\gd$ can be identified in a natural way with the whole tree $\gd$.
Let $\ell^2_r(\Gamma_x)$ denote the set of all functions in $\ell^2(\Gamma_x)$
which are radial on the subtree $\Gamma_x$.
Hence $\ell^2_r(\Gamma_x)$ can be identified with $\ell^2_r(\gd)$.

In this way the matrix $J$ restricted to $\ell^2(\Gamma_x)$ coincides with the operator $J_k$ on $\ell^2(\gd)$.
Moreover, $J$ restricted to $\ell^2_r(\Gamma_x)$ coincides with the operator $J_k$ on $\ell^2_r(\gd)$.
Similarly as in (\ref{recurrence_j_r}) and (\ref{matrix_j_r}), it can be further identified with $J_k$ on $\ell^2(\N)$.

\bigskip

From now on we make the assumption: the operator $J$ in $\ldwa$ is not essentially selfadjoint.
Hence
\[
N_z(J)\ =\ \big(\,\obraz(J-\bar{z}I)\,\big)^{\perp}
\ \neq \ \{0 \}.
\]

By Theorem \ref{iff}, the operator $J^r$ on $\ell^2(\N)$ is not essentially selfadjoint.
Furthermore, from the proof of this theorem, the truncated matrix $J_k^r$ on $\ell(\N)$
is neither essentially selfadjoint.
By the above arguments, $J$ on $\ell^2_r(\Gamma_x)$ is not essentially selfadjoint either.
Moreover, its deficiency space is one dimensional (cf. Theorem \ref{jednowymiarowy_defekt}).

Let $\tilde{f_x}$ denote a nonzero function in this deficiency space.
Observe that, $\tilde{f_x}(x)\neq 0$.
Indeed, if $v_n$ denote the value of $\tilde{f_x}$ on the $n$th level of $\gd \supset \Gamma_x$,
then the condition describing the deficiency space
\[
J\, \tilde{f_x}\ =\ z \, \tilde{f_x}
\]
(cf. Fact \ref{def}) is equivalent to the system of equations
\begin{align*}
zv_k & = \beta_k v_k + d\cdot \la_k v_{k+1},\\
zv_n & = \la_{n-1} v_{n-1} +\beta_n v_n +d\cdot \la_n v_{n+1}, \quad n>k.
\end{align*}
Hence, if $v_k=0$,  then $v_{k+1}=0$. This implies that $\tilde{f_x} \equiv 0$, which yields a contradiction.

\bigskip

Let us choose a function $\tilde{f_x}$ such that $\tilde{f_x}(x)=1$.
For each vertex $x\in\gd$ we   define the corresponding function $f_x\in\ell^2(\gd)$ by
saying that $\supp f_x \subseteq \Gamma_x$ and $f_x$ coincides with $\tilde{f_x}$ on $\Gamma_x$.

\bigskip

For each vertex $x\in\gd$ we also  define the corresponding linear subspace
\[
A_x\ = \ \bigg\{ \, \sum\limits_{i=1}^{d} a_i \cdot f_{x_i} \colon \qquad a_i\in\C \ , \qquad \sum\limits_{i=1}^{d} a_i=0 \   \bigg\}.
\]

For $i\neq j$ functions $f_{x_i}$ and $f_{x_j}$ are orthogonal as their supports are disjoint.
Note that the condition $\sum\limits_{i=1}^{d} a_i=0$ guarantees that each element  $g\in A_x$ ($x\in \gd$ and $|x|=n$)
satisfies, in addition,  the recurrence relation (\ref{rownania_defektu}) at the vertex $x$, namely
\begin{align*}
0 \ &=\ z \cdot g(x) \\
&=\  \lambda_{n-1} g(x_0) + \beta_{n} g(x)+\lambda_{n}\big( g(x_1)+g(x_2)+\ldots + g(x_d)\big) \\
 &=\ 0\ +\ 0\ +\ \la_n \,\sum\limits_{i=1}^{d} a_i\ =\ 0.
\end{align*}
It means that all  $A_x$ are $(d-1)$ dimensional  linear subspaces of the deficiency space $N_z(J)$.

\bigskip

In addition, we set
\[
A_0\ = \ \big\{\, a\cdot f_e\colon\qquad a\in\C\, \big\}.
\]
Obviously, the space $A_0$ is a one dimensional linear subspace of the deficiency space $N_z(J)$.

\bigskip

We are going to exhibit some properties of the spaces $A_x$.
First, we establish the following technical lemma.

\bigskip

\bigskip

\begin{lem}\label{znikanie_na_pietrach}
Let $x\in \gd$ and $|x|=n$. If $g\in A_x$, then
\[
\sum\limits_{\stackrel{y\in\Gamma_x}{|y|=k}} \; g (y)\ =\ 0
\]
for all $k\ge n+1$.
\end{lem}

\bigskip

\begin{dowod}
It is sufficient to make the following observation. For two different vertices $x_i$ and $x_j$ with the same predecessor $x$ the
values of the functions $f_{x_i}$ and $f_{x_j}$ on the corresponding levels of the subtrees $\Gamma_{x_i}$ and $\Gamma_{x_j}$ are equal.
This is because, by definition, the values of a function $f_y$ depend only on
the length of the vertex $y$ and on the length of the current one.
It follows that the sum of values of the function $g=f_{x_i}-f_{x_j}$ vanish on each level of the subtree $\Gamma_x$.
It is easily seen that any function $g\in A_x$ is a linear combination of functions $f_{x_i}-f_{x_j}$.
Therefore values of $g\in A_x$ also vanish on all levels of the subtree $\Gamma_x$.\\

\end{dowod}

\bigskip
\bigskip

\begin{f}\label{ortogonalnosc}
Let $x,y\in\gd\cup\{0\}$ and $x \neq y$. Then  $A_x\ \bot \ A_y$.
\end{f}

\bigskip

\begin{dowod}
Let  $g_x\in A_x$ for some vertex $x\in\gd$, where $|x|=n$.
Since $f_e$ is radial, we write $f_e(|t|)=f_e(t)$ for $t\in\gd$.
Then
\[
(\;g_x\; ,\; f_e\;)\ =\ \sum\limits_{k=n+1}^{\infty} \sum\limits_{\stackrel{t\in\Gamma_x}{|t|=k}} g_x(t) \overline{f_e(t)}
\ = \ \sum\limits_{k=n+1}^{\infty} \overline{f_e(k)} \sum\limits_{\stackrel{t\in\Gamma_x}{|t|=k}}g_x(t).
\]
By Lemma \ref{znikanie_na_pietrach}, all the sums $\sum\limits_{\stackrel{t\in\Gamma_x}{|t|=k}} g_x(t)$ vanish, whence
\[
(\;g_x\;,\;f_e\;)\ = \ 0.
\]

Consider a function $g_y\in A_y$ for some vertex $y$ different from $x$.
If  $x\notin\Gamma_y$ and $y\notin\Gamma_x$, then functions $g_x$ and $g_y$ have disjoint
supports and thus they are orthogonal.
On the other hand, if $x\in\Gamma_y$, then
\[
|x|>|y|\qquad\textrm{ and }\qquad \supp (g_x)\subset\Gamma_y.
\]
Hence
\[
(\;g_x\;,\;g_y\;)
\ =\ \sum\limits_{t\in\Gamma_x } g_x(t)\;\overline{g_y(t)}
\]
and on levels of $\Gamma_x$ the function $g_y$ has constant values $g_y(k)$.
Therefore, applying Lemma \ref{znikanie_na_pietrach} once more, we obtain
\[
(\;g_x\;,\;g_y\;) = \sum\limits_{k=n+1}^{\infty}\;\sum\limits_{\stackrel{t\in\Gamma_x}{|t|=k} } g_x(t)\overline{g_y(t)}
= \sum\limits_{k=n+1}^{\infty}\overline{g_y(k)}\sum\limits_{\stackrel{t\in\Gamma_x}{|t|=k} } g_x(t)\ =\ 0.
\]
Clearly, the case when $y\in\Gamma_x$ is similar.

\end{dowod}

\bigskip
\bigskip

\begin{f}
Assume that $f\in N_z(J)$ and $f\ \bot\ A_x$ for all $\ x\in\gd\cup\{0\}$.
Then $f\equiv 0$.
\end{f}

\bigskip

\begin{dowod}
We are going to show that such a function $f$ vanishes on the successive levels of $\gd$ starting from the root $e$.
The function $f_e$ is radial on $\gd$, whence  $E(f_e)=f_e$ (cf. (\ref{E})).
By Lemma~\ref{kontrakcja}, we thus get
\[
0\ =\ (\;f\;,\;f_e\;)\ = \ (\;f\;,\;E(f_e)\;)\ = \ (\;E(f)\;,\;f_e\;).
\]
By the same lemma,  the function $E(f)$ is square-summable.
Moreover, both $f_e$ and $E(f)$ are elements of the deficiency space $N_z(J^r)$
because taking the mean value on levels does not affect the recurrence relation (\ref{defekt_r}).
But there is a unique up to a constant function in $N_z(J^r)$.
Therefore, $E(f)$ is a constant multiple of $ f_e $.
Let $E(f)=\al f_e$. Then
\[
0\ =\ (\;E(f)\;,\;f_e\;)\ = \ (\;\alpha f_e\;,\;f_e\;)\ =\ \alpha \|f_e\|^2,
\]
whence $\alpha=0$. Thus $E(f)= 0$ and in particular
\[
f(e) \ =\ (Ef)(e)\ =\ 0.
\]
The orthogonality of $f$ to the function $f_e$ yields that $f$ vanishes at the root $e$,
i.e. on the zero level of the tree $\gd$.
The orthogonality of $f$ to the successive spaces $A_x$ enables us to show that $f$ is
equal to zero at the corresponding vertices.
Assume that $f(x)=0$ for each $|x|\le n$.
Fix a vertex $x$ of length  $n$.
Since $f\in N_z(J)$ and $f(x)=f(x_0)=0$, the recurrence equation (\ref{rownania_defektu}) at $x$
\[
z f(x) = \la_{n-1}f(x_0)  +  \beta_nf(x) + \la_n (f(x_1)+f(x_2)+\ldots+f(x_d))
\]
gives 
\begin{equation}\label{a}
f(x_1)+f(x_2)+\ldots+f(x_d)\ =\ 0.
\end{equation}
Fix $g\in A_x$.
Since the function $g$ is radial on each subtree $\Gamma_{x_i}$,
taking the mean value on each one separately (cf. (\ref{E_x})) does not affect the values of $g$, i.e.
\[
E_{x_d} \;E_{x_{d-1}} \;\ldots \; E_{x_1}\;( g)\ = \ g.
\]
By the symmetry of these averaging operators (cf. Lemma \ref{kontrakcja}),
\begin{align*}
0\ = \ (\;f\;,\;g\;)\
&= \ (\;f\;,\; E_{x_d} E_{x_{d-1}} \ldots  E_{x_1}\;g\;) \\
&= \ (\;E_{x_1} E_{x_{2}} \ldots  E_{x_d}\;f\;,\;g\;).
\end{align*}
By (\ref{a}),  the function
\[
\mathds{1}_{\Gamma_{x}}\cdot E_{x_1} E_{x_{2}} \ldots  E_{x_d}\;f,
\]
where  $\mathds{1}_{\Gamma_{x}}$ denotes the characteristic function of $\Gamma_x \supseteq \supp g$,
belongs to and is orthogonal to the space $A_x$ at the same time. Hence it must be zero.
Therefore
\[
f(x_i)\ =\ \mathds{1}_{\Gamma_{x}}\cdot E_{x_1}E_{x_2}\ldots E_{x_d}f(x_i)\ = \ 0
\]
\nopagebreak
for all $i=1,2,\ldots,d$.
\end{dowod}

\bigskip
\bigskip

We thus see that the sets $A_x$, in a sense, fill up the whole of the deficiency space  $N_z(J)$.
To be more precise, the above facts can be summarized as follows.

\bigskip
\bigskip

\begin{tw}\label{suma_a_x}
The algebraic direct sum
\[
\bigoplus\limits_{x\in\gd\cup\{0\}} A_x \ =\
\lin \;\big\{ \; g_x\colon\quad g_x\in A_x, \ x\in\gd\cup\{0\} \big\}
\]
of the pairwise orthogonal spaces $A_x$
is dense in the nontrivial deficiency space $N_z(J)$.
\end{tw}

\bigskip
\bigskip

\noindent
{\textbf{Remark.}}
In the case when $d=2$ not only $A_0$ but also all the remaining sets $A_x$ for $x\in\Gamma_2$
are one dimensional linear spaces. Moreover, the functions
\[
g_x \ =\ \dfrac{f_{x_1}-f_{x_2}}{\|f_{x_1}-f_{x_2}\|}\;, \qquad\qquad x\in\gd,    
\]
along with the function $g_0=\frac{f_0}{\|f_0\|}$
form the orthonormal basis in $N_z(J)$ on the tree $\Gamma_2$.

\bigskip
\bigskip


Let us now calculate norms of elements of $A_x$ in the case  when  $d\ge 2$ is arbitrary.

\bigskip

Let $p_n$ be the orthogonal polynomials (cf. (\ref{a_n}))  associated with the matrix
\begin{equation}\label{p_n}
J^r\ =\
\left(
\begin{array}{cccccc}
\beta_0   & \pd\la_0 & 0 & 0 & 0  & \ldots\\
\pd\la_0  & \beta_1   & \pd\la_1 & 0 & 0 & \ldots\\
0 & \pd\la_1  & \beta_2   & \pd\la_2 &  0 & \ldots\\
0 & 0 & \pd\la_2  & \beta_3   & \pd\la_3 &   \ldots\\
0&0&0& \pd\la_3 & \beta_4 & \ddots \\
\vdots & \vdots & \vdots & \vdots & \ddots & \ddots
\end{array}
\right),
\end{equation}
i.e. let the numbers $p_n(z)$ satisfy the equations
\begin{equation}\label{gwiazdka}
zp_n(z)=\sqrt{d}\la_{n-1}p_{n-1}(z)+\beta_np_n(z)+\sqrt{d}\la_np_{n+1}(z),\quad n\ge 0,
\end{equation}
\[
p_{-1}(z)\ =\ 0,\qquad p_0(z)\ =\ 1.
\]
Dividing by $\big(\sqrt{d}\big)^{n}$ gives
\[
z \cdot \dfrac{p_n(z)}{\sqrt{d^n}}= \la_{n-1}\cdot \dfrac{p_{n-1}(z)}{\sqrt{d^{n-1}}}
+\beta_n\cdot\dfrac{p_n(z)}{\sqrt{d^n}}+d\cdot\la_n\cdot\dfrac{p_{n+1}(z)}{\sqrt{d^{n+1}}} ,\quad n\ge 0.
\]
These are exactly the equations describing the unique radial function in $N_z(J)$ (cf. (\ref{defekt_r})), hence
\begin{equation}\label{wartosci_0}
f_0(x)\ =\ f_0(|x|)\ =\ \dfrac{p_{|x|}(z)}{\sqrt{d^{|x|}}}.
\end{equation}
We calculate the norm
\[
\|f_0\|^2\ =\ \sum\limits_{n=0}^{\infty} \; d^n\; \big|f_0(n)\big|^2\ =\ \sum\limits_{n=0}^{\infty} \; d^n\;  \bigg|\dfrac{p_n(z)}{\sqrt{d^n}}\bigg|^2
\ =\ \sum\limits_{n=0}^{\infty} \;|p_n(z)|^2.
\]
The norm  of an arbitrary function $g\in A_0$ can be expressed as
\[
\|g\|\ =\ \alpha_{0}(z) \cdot |g(e)|,
\]
where
\[
\alpha_{0}(z)\ =\ \bigg( \, \sum\limits_{n=0}^{\infty} \;|p_n(z)|^2\,\bigg)^{\frac{1}{2}}.
\]

\bigskip

Let $q_n$ be the orthogonal polynomials of the second kind (cf. (\ref{a_n})) associated with the matrix $J^r$, i.e.
\[
zq_n(z)=\sqrt{d}\la_{n-1}q_{n-1}(z)+\beta_nq_n(z)+\sqrt{d}\la_nq_{n+1}(z),\quad n\ge 1,
\]
\[
 q_0(z)\ =\ 0, \qquad q_1(z)=\frac{1}{\la_0}.
\]
As before, dividing by $\big(\sqrt{d}\big)^{n-1}$ gives
\[
z \cdot \dfrac{\la_0 q_n(z)}{\sqrt{d^{n-1}}}= \la_{n-1}\cdot \dfrac{\la_0 q_{n-1}(z)}{\sqrt{d^{n-2}}}
+\beta_n\cdot\dfrac{\la_0 q_n(z)}{\sqrt{d^{n-1}}}+d\cdot\la_n\cdot\dfrac{\la_0 q_{n+1}(z)}{\sqrt{d^{n}}} ,\quad n\ge 1.
\]
Therefore, for a fixed $i$, the  values of $f_{e_i}$ are expressed by polynomials $q_n$ as follows
\begin{equation}\label{wartosci_1}
f_{e_i}(x)\ =\ f_{e_i}(|x|)\ =\ \la_0\cdot\dfrac{q_{|x|}(z)}{\sqrt{d^{|x|-1}}}.
\end{equation}
Hence
\begin{align*}
\|f_{e_i}\|^2
&=\sum\limits_{x\in\Gamma_{e_i}} |f_{e_i}(x)|^2
= \sum\limits_{n=1}^{\infty}\; d^{n-1}\; \big|f_{e_i}(n)\big|^2\\
&= \sum\limits_{n=1}^{\infty} \;d^{n-1}\;\la_0^2 \; \bigg|\dfrac{q_n(z)}{\sqrt{d^{n-1}}}\bigg|^2
 = \la_0^2 \;\sum\limits_{n=1}^{\infty}\; |q_n(z)|^2.
\end{align*}
Let 
\begin{equation}\label{alfa}
\alpha_1(z)\ =\ \la_0 \bigg(\,\sum\limits_{n=1}^{\infty}\; |q_n(z)|^2\,\bigg)^{\frac{1}{2}}.
\end{equation}
Since the functions $f_{e_i}$ are pairwise orthogonal,
the norm of an arbitrary function $g\in A_e$ is equal to
\[
\|g\|\ =\  \Big\| \sum\limits_{i=1}^d g(e_i)f_{e_i} \Big\| \ =\ \alpha_1(z)\cdot\sqrt{|g(e_1)|^2+\ldots+|g(e_d)|^2}.
\]

\bigskip

Now we consider a function $g\in A_x$ for a fixed vertex $x\neq e$, i.e. $|x|=k \ge1$.
Since
\[
g\ = \ \sum\limits_{i=1}^d g(x_i) \cdot f_{x_i},
\]
where $f_{x_i}$ are pairwise orthogonal, we get
\[
\| g \|^2 \ = \ \|f_{x_i}\|^2\cdot \sum\limits_{i=1}^d |g(x_i)|^2,
\]
because the values of $f_{x_i}$ on the subtree $\Gamma_{x_i}$ depend only on the length of vertices
and on $k$ which is the length of the root of this subtree.
These values are determined by the equations
\[
z f_{x_i}(n)=\la_{n-1}f_{x_i}(n-1)+\beta_n f_{x_i}(n)+ d \la_n f_{x_i}(n+1),\qquad n\ge k+1,
\]
and
\[
f_{x_i}(k)=0\ ,\qquad\qquad f_{x_i}(k+1)=1
\]
(cf. (\ref{defekt_r}) and the definition of $f_x$).
Note that the numbers
\[
\frac{\la_k\; \big(\;p_k(z)q_n(z)-q_k(z)p_n(z)\;\big)}{\sqrt{d^{n-(k+1)}}},\qquad\quad n\ge k,
\]
satisfy these equations. Indeed,
the recurrence relation results from the fact that the orthogonal polynomials of the first and  second kind, i.e. numbers $p_n(z)$ and $q_n(z)$,
satisfy the recurrence starting with $n=1$, in particular, for $n\ge k$. Therefore, the same holds for any linear combination of them.
Since $|x_i|=k+1$,  there are exactly $n-(k+1)$ vertices on the $n$th level in the subtree $\Gamma_{x_i}$.
This is why the exponent of the power in the denominator is equal to $n-(k+1)$.
Furthermore, for $n=k$ the value is 0.
Finally, by the  formula
\[
 p_n(z)q_{n+1}(z)-p_{n+1}(z)q_n(z)\ =\ \frac{1}{\la_n}
\]
connecting the polynomials of the first and  second kind (see e.g. \cite{akhiezer} or \cite{szwarc}),
we get the value 1 for $n=k+1$.
Consequently,
\begin{equation}\label{wartosci}
f_{x_i}(n)\ = \ \frac{\;\la_k \big(\;p_k(z)q_n(z)-q_k(z)p_n(z)\;\big)\;}{\sqrt{d^{n-(k+1)}}}, \qquad n\ge k+1,
\end{equation}
and thus
\begin{align*}
\|f_{x_i}\|^2 \
 &= \ \sum\limits_{n=k+1}^{\infty} d^{n-(k+1)}\;|f_{x_i}(n)|^2\\
 & =\ \la_k^2\sum\limits_{n=k+1}^{\infty}\big|p_k(z)q_n(z)-q_k(z)p_n(z)\big|^2.
\end{align*}
Hence
\[
\|g\|^2 =  \Big(|g(x_1)|^2+\ldots+|g(x_d)|^2 \Big)\cdot\la_k^2\sum\limits_{n=k+1}^{\infty}\big|p_k(z)q_n(z)-q_k(z)p_n(z)\big|^2.
\]
Let $\alpha_{k+1}(z)$ denote the positive number such that
\begin{equation}\label{alfy}
\alpha_{k+1}^{\ 2}(z)\ =\ \la_k^2\sum\limits_{n=k+1}^{\infty}\big|p_k(z)q_n(z)-q_k(z)p_n(z)\big|^2, \qquad k\ge1.
\end{equation}
Then
\[
\|g\|\ = \ \alpha_{k+1}(z) \cdot\sqrt{|g(x_1)|^2+\ldots+|g(x_d)|^2}
\]
for any function $g\in A_x$, where $|x|=k\ge 1$.

\bigskip

Note that  for $k=0$ the right hand side of (\ref{alfy}) gives exactly
the number $\alpha_1(z)$ defined already  by  (\ref{alfa}) so the
numbers $\alpha_k(z)$ may be defined by  the common formula (\ref{alfy}) for all $k\ge 0$.

\bigskip

The following fact is a summary of the previous considerations concerning norms.

\bigskip

\begin{f} \label{normy}
We have
\[
\|f_x\|\ = \ \al_{|x|} \qquad\quad \textrm{ for } \qquad x\in\gd \cup \{0\},
\]
and
\[
\|g\|\ =\
\left\{
\begin{array}{ll}
\alpha_0(z)\ |g(e)|\ ,&\ \textrm{ if }\  g\in A_0,\\ & \\
\alpha_{|x|+1}(z)\ \Big( \sum\limits_{i=1}^{d} |g(x_i)|^2 \Big)^{\frac{1}{2}}\ ,&\ \textrm{ if }\ g\in A_x,\ x\in\gd,
\end{array}\right.
\]
where the coefficients $\alpha_k(z)$ do not depend on functions and are expressed as follows
\begin{align*}
\alpha_0(z)\ &=  \sum\limits_{n=0}^{\infty} \;|p_n(z)|^2 ,\\
\alpha_k(z)\ &= \ \la_{k-1}^2\sum\limits_{n=k}^{\infty}\big|p_{k-1}(z)q_n(z)-q_{k-1}(z)p_n(z)\big|^2  \quad\textrm{ for } \ k\ge 1.
\end{align*}
\end{f}

\bigskip



\clearpage

\hspace{7pt}

\begin{center}
{\textsc{The deficiency space and the boundary of tree}}
\end{center}

\hspace{7pt}

A \emph{path} in a tree is, by definition,
a sequence $\{x_n\}$ of vertices such that for any $n$, the vertices $x_n$ and $x_{n+1}$
are joined by an edge.
The \emph{boundary} $\Omega=\bgd$ of the tree $\gd$ is the set of all infinite paths starting at the root $e$.

Note that at each level
on the way downward from the root $e$  we have to choose one of $d$ edges,
hence the boundary $\Omega$ can be identified with the Cantor set
\[
\Omega \ \simeq  \{0,1,2,\ldots,d-1\}^{\mathbb{N}}
\]
(which is the classical Cantor set in an interval when $d=2$).
Clearly, each vertex $x$, and thereby each subtree  $\Gamma_x$, is associated with
a \emph{cylindric set (a cylinder)} $\Omega_x\subseteq \Omega$, i.e.
the set of all those paths which contain the vertex~$x$.

Let $\mu$ be the probability measure on $\{0,1,2,\ldots, d-1\}$ such that
\[
\mu \ = \ \frac{1}{d}\cdot\bigg(\delta_0 + \delta_1+\ldots+\delta_{d-1}\bigg).
\]
Let $\dw$ denote the natural probability product measure on the boundary $\Omega$
\[
\dw\ =\ \bigotimes\limits_{i=0}^{\infty} d\mu_i,\qquad \qquad \mu_i=\mu,
\]
i.e. values of $\dw$ on cylindric sets are given by
\[
\dw(\;\Omega_x\;) \ = \ d^{-|x|},\qquad
\qquad x\in\gd.
\]

\bigskip

We consider the space $L^2(\Omega, \dw)$ of those functions defined on the boundary $\Omega$
which are square-summable with respect to the measure~$\dw$.
For each subspace $A_x\subset N_z(J)$ we define the corresponding subspace
$B_x \subseteq L^2(\Omega,\dw)$.
Namely, let $B_0$ denote the one dimensional linear subspace of constant functions on $\Omega$
and for $x\in\gd$ we put
\[
B_x\ =\ \bigg\{\; \sum\limits_{i=1}^d b_i\cdot\1{\Omega_{x_i}}\colon\quad b_i\in\C, \quad\sum\limits_{i=1}^{d} b_i = 0\ \bigg\}.
\]
Similarly to $A_x$, each subspace $B_x$ is a linear space of dimension $d-1$.
The another analogy is given by the following property of any element $F$ of the space $B_x\,$:
\[
\int\limits_{\Omega} F(\omega) \dw \ = \ \sum\limits_{i=1}^d \int\limits_{\Omega_{x_i}} F(\omega) \dw \
= \ d^{-|x_i|} \cdot \sum\limits_{i=1}^d b_i \ = 0.
\]

\bigskip
\bigskip

\begin{f}\label{B_x}
The subspaces $B_x$ for $x\in\gd\cup\{0\}$  are pairwise orthogonal and fill up the whole of the space
$L^2(\Omega, \dw)$, i.e. the algebraic direct sum
\[
\bigoplus\limits_{x\in\gd\cup\{0\}} B_x \ =\
\lin \;\big\{ \; G_x\colon\quad G_x\in B_x, \ x\in\gd\cup\{0\} \big\}
\]
is dense in $\Ldwa$.
\end{f}

\bigskip

\begin{dowod}
Let $G_x\in B_x$ and $G_y\in B_y$ for $x\neq y$.
Cylindric sets associated with two different vertices are either disjoint or one is a proper subset of another.
If $\Omega_x \cap \Omega_y = \emptyset$, then functions $G_x$ and $G_y$ are orthogonal as their supports are disjoint.
On the other hand, if $\Omega_y \varsubsetneq \Omega_x$, then there exists $i$ such that $\Omega_y \subseteq \Omega_{x_i}$.
Let $b_i$ denote the value of $G_x$ on the cylinder $\Omega_{x_i}$. Then
\[
\iloczyn{G_x}{G_y}\ =\ \int\limits_{\Omega_{y}}G_x(\omega)G_y(\omega)\dw
= \ b_i \cdot \int\limits_{\Omega_y} G_y(\omega) \dw\ = \ 0,
\]
which completes the proof of orthogonality.

Assume that a function $F\in \Ldwa$ is orthogonal to every $B_x$ for $x\in\gd\cup \{0\}$.
In particular, for the function $G_0 \equiv 1$ belonging to the space $B_0$ we obtain
\begin{equation}\label{calka_zero}
0\ =  \ \iloczyn{F}{G_0} \ = \ \int\limits_{\Omega} F(\omega) \dw\ =\ \sum\limits_{i=1}^d \;\,\int\limits_{\Omega_{e_i}} F(\omega) \dw.
\end{equation}
The orthogonality of $F$ to $\1{\Omega_{e_i}}-\1{\Omega_{e_j}}\in B_e$ for $i\neq j$ gives
\[
0\ = \ \iloczyn{F}{\1{\Omega_{e_i}}-\1{\Omega_{e_j}}}\ = \ \int\limits_{\Omega_{e_i}} F(\omega) \dw - \int\limits_{\Omega_{e_j}} F(\omega) \dw,
\]
whence
\[
\int\limits_{\Omega_{e_i}} F(\omega) \dw\ =\ \int\limits_{\Omega_{e_j}} F(\omega) \dw.
\]
Since all the numbers
\[
\int\limits_{\Omega_{e_i}} F(\omega) \dw
\]
are equal and sum up to 0 (the equality (\ref{calka_zero})), all of them vanish.
Similar considerations applied to $x=e_j$ and its successors $x_i$ yield
\[
\int\limits_{\Omega_y} F(\omega)\dw \ = \ 0
\]
for $\dist (y, e) =2$.
In this way one can show  that integrating $F$ over an arbitrary cylindric set gives 0.
Hence $F=0$\ \ $\dw$--almost everywhere.

\hspace{10pt}
\end{dowod}

\bigskip
\bigskip

Let $x\in\gd$.
For the function $f_x$ we define the corresponding function $F_x\in\Ldwa$ by
\[
F_x\ = \ \al_{|x|} \cdot\sqrt{d^{|x|}}\cdot \1{\Omega_{x}}.
\]
Therefore
\begin{align*}
\|F_x\|^2\
= \ \int\limits_{\Omega} F_x(\omega)\overline{F_x(\omega)} \dw
\ &= \ \al^2_{|x|}(z) d^{|x|} \cdot \dw(\Omega_x)\\
&= \ \al^2_{|x|}(z)\ = \ \|f_x\|^2.
\end{align*}
Clearly, if $x_i\neq x_j$ have a common predecessor, then functions
$F_{x_i}$ and $F_{x_j}$ are orthogonal as their supports are disjoint.
Thus the formula
\begin{equation}\label{odpowiednie}
A_x\ \ni\ \sum\limits_{i=1}^d a_i f_{x_i}\ =\ g \ \ \ \longleftrightarrow\ \ \ G\ =\ \sum\limits_{i=1}^d a_i F_{x_i}\ \in\  B_x
\end{equation}
sets a one-to-one correspondence between functions in $A_x$ and functions in $B_x$.
Furthermore,
\begin{align*}
\| G\|\ &=\ \| F_{x_i}\| \cdot \sqrt{|a_1|^2+|a_2|^2+\ldots+|a_d|^2 } \\
        &=\ \| f_{x_i}\| \cdot \sqrt{|a_1|^2+|a_2|^2+\ldots+|a_d|^2 }\ =\ \|g\|.
\end{align*}
It follows that a mapping
\[
 G\ \longmapsto\  g
\]
is a linear bijection between the spaces $B_x$ and $A_x$ which, in addition, preserves the norm.
Apart from these bijections for $x\in\gd$, we set the mapping from  $A_0$  onto $B_0$ by
\begin{equation}\label{odpowiednie2}
f_0 \quad \longmapsto \quad F_0=\al_0\cdot \1{\Omega},
\end{equation}
which, clearly, is also a norm preserving linear bijection.

\bigskip

\noindent
In view of Theorem \ref{suma_a_x} and Fact \ref{B_x}, all these bijections have a unique extension
to the injective isometry
\begin{equation}\label{izometria}
U\colon\quad L^2(\Omega,\dw) \ \ \stackrel{\textrm{onto}}{\longrightarrow}\ \ N_z(J).
\end{equation}

\bigskip
\bigskip

For a fixed vertex $y\in\gd$ we define the functional on $L^2(\Omega,\dw) $ by
\[
F\ \longmapsto\ \big(UF\big)(y).
\]
As it is linear and  bounded, it determines, by the Riesz Theorem, a unique function  $P_z(y,\omega)\in L^2(\Omega,\dw)$
such that for all functions $F\in L^2(\Omega,\dw)$
\begin{equation}\label{riesz}
\big(UF\big)(y)\ = \ \int\limits_{\Omega} P_z(y,\omega) F(\omega) \dw.
\end{equation}
In particular, for a function $g\in A_x$ and the corresponding function $G\in B_x$ we have
\begin{equation}\label{poisson}
g(y)\ =\  \int\limits_{\Omega} P_z(y,\omega) G(\omega) \dw.
\end{equation}
In view of this formula, it is natural to call the function $P_z(y,\omega)$ the
{\emph{Poisson kernel}}.
It describes a relationship between functions in the deficiency space and functions on the boundary of the tree.

\bigskip

We turn to the description of the Poisson kernel $P_z(y,\omega)$.
In order to do this we state some of the properties of this kernel.

\bigskip
\bigskip

\begin{f}
For a fixed $y\in\gd$
\[
\big( JP_z(\cdot , \omega) \big)(y)\ = \ z\cdot P_z(y, \omega)
\]
 $\dw$--almost everywhere.
\end{f}

\bigskip

\begin{dowod}
Let $y\in\gd\cup\{0\}$.
It is sufficient to show that for an arbitrary $x\in\gd\cup\{0\}$ and an arbitrary function $G\in B_x$
\[
(\, JP_z(y, \cdot)\, , \, G\, )\ =\ (\, zP_z(y,\cdot)\, , \, G \,).
\]
Let us consider the function $g = UG \in A_x$.
Since $A_x\subset N_z(J)$,  $g$ satisfies, in particular, the recurrence relation (\ref{rownania_defektu}) at $y$, i.e.
\[
(J g) (y) \ =\ z \cdot g(y).
\]
In view of (\ref{poisson}), we have
\[
\bigg( J\Big( \int\limits_{\Omega}P_z(\cdot,\omega)G(\omega)\dw\Big)\bigg)(y)\
= \ z \cdot \int\limits_{\Omega}P_z(y,\omega)G(\omega)\dw.
\]
By the linearity of the integral, we get
\[
\int\limits_{\Omega}JP_z(y,\omega)G(\omega)\dw = \int\limits_{\Omega}zP_z(y,\omega)G(\omega)\dw,
\]
hence
\begin{align*}
\big(\, JP_z(y, \cdot)\, , \, \overline{G}\, \big)\ =\ \big(\, zP_z(y,\cdot)\, , \, \overline{G} \,\big).
\end{align*}
\end{dowod}

\bigskip
\bigskip

Now we are about to use automorphisms of the tree $\gd$.

\bigskip
\bigskip

\noindent
{\textbf{Remark.}}
Each automorphism of the tree $\gd$ leaves the root $e$ fixed.
Indeed, since any automorphism maps any vertex to a vertex of the same degree
(i.e. the number of edges to which a given
vertex  belongs), it suffices to observe that
the root $e$ is the only one vertex in $\gd$ of degree $d$.
Obviously, each automorphism of $\gd$ acts on the boundary $\Omega$ at the same time.

\bigskip
\bigskip

Note that $J$ commutes with all isometries of the tree $\gd$.
Each automorphism of $\gd$ acts on functions in $A_x$ and in $B_x$ in a natural way.
Namely, let $k\colon\gd\to\gd$ be an automorphism, i.e. $k$ satisfy
\[
(\;\forall x,y\in\gd\;)\  \ \ \dist(kx,ky) \ = \ \dist(x,y),
\]
\[
ke\ =\ e.
\]
For any $g\in A_x$ we put
\[
(_kg) (y) \ =\ g( k^{-1}y).
\]
Then
\[
_kg\in A_{kx}\qquad\ \textrm{ and }\ \qquad _kg(ky)=g(y) .
\]
Similarly, $k$ acts on functions in $B_x$. Hence for the corresponding
function $G\in B_x$ we have
\[
_kG\in B_{kx}\qquad\ \textrm{ and }\ \qquad (_kG)(k\omega)\ = \ G(\omega).
\]

\bigskip
\bigskip

\begin{lem}\label{izometria_a_jadro}
For a fixed vertex $y\in\gd$ and an arbitrary automorphism $k$ of the tree $\gd$
\[
P_z(\, ky\,,\, k\omega\,) \ = \  P_z(\, y\,,\, \omega\,)
\]
 $\dw$--almost everywhere.
\end{lem}

\bigskip

\begin{dowod}
Let $x\in\gd\cup\{0\}$ and $g\in A_x$.
By (\ref{poisson})  and the property of the automorphism $k$, we have
\begin{align*}
\int\limits_{\Omega} P_z(y,\omega) G(\omega) \dw
 & = g(y)
  = {_k}g (ky) = U(_kG)\,(ky) \\
 &= \int\limits_{\Omega} P_z(ky,\omega)\cdot (_kG)(\omega) \,\dw .
\end{align*}
Replacing  $\omega$ by $k\omega$ in the last integral yields
\begin{align*}
\int\limits_{\Omega} P_z(y,\omega)\cdot G(\omega) \dw
& =  \int\limits_{\Omega} P_z(ky,k\omega)\cdot (_kG)(k\omega) \;d(k\omega)\\
& =  \int\limits_{\Omega} P_z(ky,k\omega)\cdot G(\omega) \;d(k\omega).
\end{align*}
By the invariance of the measure $\dw$, we obtain
\begin{align*}
\int\limits_{\Omega} P_z(y,\omega)\cdot G(\omega) \dw
& =  \int\limits_{\Omega} P_z(ky,k\omega)\cdot G(\omega) \,\dw.
\end{align*}
Since $G$ was an arbitrary function, we get
\[
P_z(ky,k\omega)=P_z(y,\omega)\qquad \dw-\textrm{almost everywhere}.
\]
\end{dowod}

\bigskip
\bigskip

Fix a path $\omega\in\Omega$ and a vertex $y\in\gd$.
Let us number the adjacent vertices in $\omega$ with the numbers $0,1,2,\ldots$ starting with the root $e$
\[
\omega\ =\ \{ \; \omega_0\,,\,\omega_1\,,\, \omega_2\,,\, \ldots \;  \},
\]
\[
\omega_0=e.
\]
The relative position of the path $\omega$ and vertex $y$ can be described by two nonnegative integers.
Let $n=n(y,\omega)$ denote the distance between $y$ and $\omega$ and let $m=m(y,\omega)$
be such that the vertex $\omega_m\in\omega$ realizes this distance $n$, i.e.
\[
n\ =\ \dist(y,\omega)\ =\ \dist(y,\omega_m).
\]
Obviously, $|y|=m+n$.
One can say that on the way from the root $e$ to the vertex $y$ we do exactly $m$ steps along the path $\omega$
and exactly $n$ steps off $\omega$.
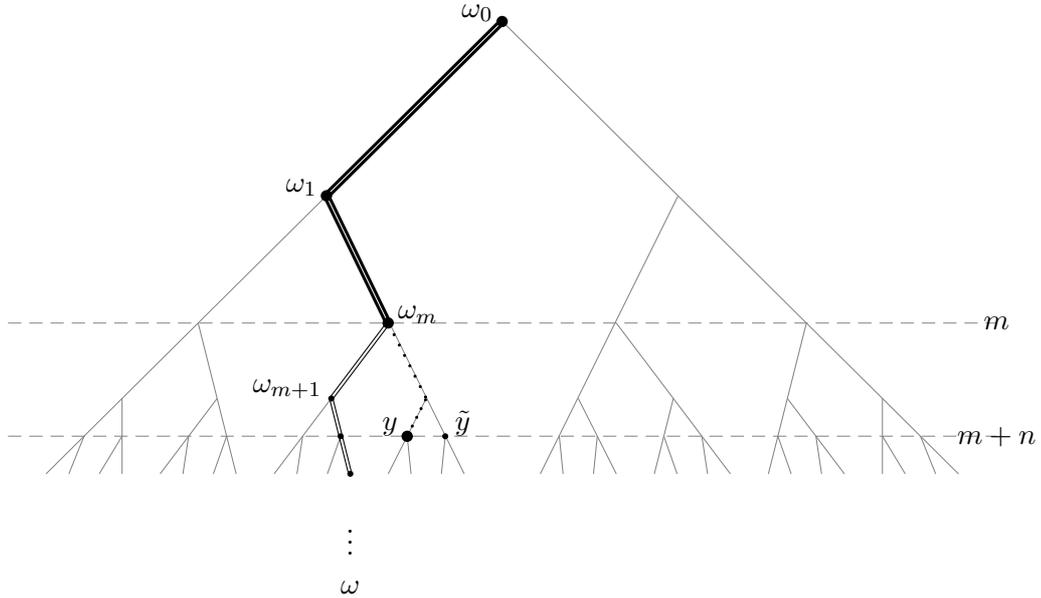
\begin{center}
\small
\begin{pspicture}(-6.3,-2)(6.3,6.5)
\psline[linewidth=0.2pt,linecolor=gray](6,0)(0,6)
\psline[linewidth=0.2pt,linecolor=gray](2.3125,3.6875)(0.5,0)

\psline[linewidth=0.2pt,linecolor=gray](1.5,2 )(3,0 )
\psline[linewidth=0.2pt,linecolor=gray](4, 2)(3.5,0 )

\psline[linewidth=0.2pt,linecolor=gray](1,1 )(1.5,0 )
\psline[linewidth=0.2pt,linecolor=gray](2.25,1 )(2,0 )
\psline[linewidth=0.2pt,linecolor=gray](3.75,1 )(4.5,0 )
\psline[linewidth=0.2pt,linecolor=gray](5,1 )(5,0 )

\psline[linewidth=0.2pt,linecolor=gray](0.75 ,0.5)(0.8 ,0)
\psline[linewidth=0.2pt,linecolor=gray](1.25 ,0.5)(1.2 ,0)
\psline[linewidth=0.2pt,linecolor=gray](2.125 ,0.5)(2.3 ,0)
\psline[linewidth=0.2pt,linecolor=gray](2.625 ,0.5)( 2.7,0)
\psline[linewidth=0.2pt,linecolor=gray](3.625 ,0.5)(3.8 ,0)
\psline[linewidth=0.2pt,linecolor=gray](4.125 ,0.5)( 4.2,0)
\psline[linewidth=0.2pt,linecolor=gray](5 ,0.5)( 5.3,0)
\psline[linewidth=0.2pt,linecolor=gray](5.5 ,0.5)(5.7 ,0)

\psline[linewidth=0.2pt,linecolor=gray](-6,0)(0,6)
\psline[linewidth=0.2pt,linecolor=gray](-2.3125,3.6875)(-0.5,0)

\psline[linewidth=0.2pt,linecolor=gray](-1.5,2 )(-3,0 )
\psline[linewidth=0.2pt,linecolor=gray](-4, 2)(-3.5,0 )

\psline[linewidth=0.2pt,linecolor=gray](-1,1 )(-1.5,0 )
\psline[linewidth=0.2pt,linecolor=gray](-2.25,1 )(-2,0 )
\psline[linewidth=0.2pt,linecolor=gray](-3.75,1 )(-4.5,0 )
\psline[linewidth=0.2pt,linecolor=gray](-5,1 )(-5,0 )

\psline[linewidth=0.2pt,linecolor=gray](-0.75 ,0.5)(-0.8 ,0)
\psline[linewidth=0.2pt,linecolor=gray](-1.25 ,0.5)(-1.2 ,0)
\psline[linewidth=0.2pt,linecolor=gray](-2.125 ,0.5)(-2.3 ,0)
\psline[linewidth=0.2pt,linecolor=gray](-2.625 ,0.5)( -2.7,0)
\psline[linewidth=0.2pt,linecolor=gray](-3.625 ,0.5)(-3.8 ,0)
\psline[linewidth=0.2pt,linecolor=gray](-4.125 ,0.5)( -4.2,0)
\psline[linewidth=0.2pt,linecolor=gray](-5 ,0.5)( -5.3,0)
\psline[linewidth=0.2pt,linecolor=gray](-5.5 ,0.5)(-5.7 ,0)

\psline[linewidth=0.2pt,linecolor=gray,linestyle=dashed](-6.5,2 )(6.5,2 )

\psline[linewidth=0.2pt,linecolor=gray,linestyle=dashed](-6.5,0.5 )(6.5,0.5 )

\psline[linewidth=1.2pt,doubleline=true,doublesep=0.3pt](0,6)(-2.3125,3.6875)
\psline[linewidth=1.2pt,doubleline=true,doublesep=0.3pt](-1.5,2)(-2.3125,3.6875)
\psline[linewidth=0.4pt,doubleline=true,doublesep=0.8pt](-2.25,1)(-1.5,2)
\psline[linewidth=0.4pt,doubleline=true,doublesep=0.8pt](-2.25,1)(-2,0)

\psline[linewidth=1.2pt,linestyle=dotted](-1,1)(-1.5,2)
\psline[linewidth=1.2pt,linestyle=dotted](-1,1)(-1.25,0.5)

\psdots[dotsize=0.15](0,6) \rput[br](0,6){$\omega_0\ $}
\psdots[dotsize=0.15](-2.3125,3.6875) \rput[br](-2.3125,3.6875){$\omega_1\ $}
\psdots[dotsize=0.15](-1.5,2) \rput[bl](-1.5,2){$\ \omega_m $}
\psdots[dotsize=0.15](-1.25 ,0.5)\rput[br](-1.25 ,0.5){$y\ $}
\psdots[dotsize=0.08](-2.25,1 )\rput[br](-2.25,1 ){$\omega_{m+1}\ $}
\psdots[dotsize=0.08](-0.75 ,0.5)\rput[bl](-0.75 ,0.5){$\ \tilde{y} $}
\psdots[dotsize=0.08](-2.125 ,0.5)(-2,0)

\rput{*0}(-2,-0.8){$\vdots$}
\rput{*0}(-2,-1.5){$\omega$}

\rput{*0}(6.5,2 ){$m$}
\rput{*0}(6.5,0.5 ){$m+n$}

\end{pspicture}
\end{center}
For instance, for $y$ and $\omega$ marked in the picture above ($d=2$)
we get $n=m=2$. Note that replacing $y$ by $\tilde{y}$ gives the same numbers $n$ and $m$.

\bigskip
\bigskip

\begin{f}
The values of the Poisson kernel $P_z(y,\omega)$ depend only on the numbers $m(y,\omega)$ and $n(y,\omega)$ defined above.
\end{f}

\bigskip

\begin{dowod}
Let $y,y'\in\gd$ and $\omega,\omega'\in\Omega$ satisfy
\[
m(y,\omega)=m(y',\omega')
\qquad \textrm{ and } \qquad
n(y,\omega)=n(y',\omega').
\]
Let $k_1$ be any automorphism on $\gd$ mapping $\omega$ to $\omega'$, i.e.
\[
k_1\,\omega \ = \ \omega'.
\]
By assumption, $\dist(k_1 y, \omega')=\dist(y,\omega')$
and both distances are realized by the same point in the path $\omega'$.
Hence, there exists an automorphism $k_2$ which fixes $\omega'$
but maps the vertices $k_1 y $ to $y'$, i.e.
\[
k_2\,\omega'\ =\ \omega' \qquad\textrm{ and }\qquad k_2(k_1 y)\ =\ y'.
\]
By Lemma \ref{izometria_a_jadro}, we thus get
\begin{align*}
P_z(y',\omega')\ &=\ P_z(k_2(k_1y), k_2\omega')\ =\ P_z(k_1 y, \omega')\\
&=\ P_z(k_1y,k_1\omega)\ =\ P_z(y,\omega).
\end{align*}
\end{dowod}

\bigskip
\bigskip

In the next part we intend to give an explicit formula for the Poisson kernel $P_z(y,\omega)$.
In order to do this we introduce some projections.

\bigskip
\bigskip

\noindent
Let $\pi_z$ denote the projection of $\ldwa$ onto the deficiency space
\[
N_z(J)=\overline{\bigoplus\limits_{x\in\gd\cup\{0\}}A_x}
\]
(cf. Theorem \ref{suma_a_x}). Then for an arbitrary function $f\in\ldwa$ we have
\[
\pi_z (f) \ = \ \sum\limits_{x\in\gd\cup\{0\}}  \pi_{z,x}(f),
\]
where $\pi_{z,x}$ denotes the projection of  $\ldwa$ onto $A_x$;
in particular $\pi_{z,0}$ is the projection onto $A_0$.

\bigskip
\bigskip

\begin{f}\label{rzuty}
Let $x$ be a vertex in $\gd$ and $k=|x|+1$. If $y\in\Gamma_x$ is different from $x$ and the number $i$ is such that
$y\in\Gamma_{x_i}$, then
\[
\pi_{z,x}(\delta_y)\ = \ \frac{\ \ \overline{f_{x_i}(y)}\ \ }{\al_k^{2}(z)} \,\cdot\, \Big[ f_{x_i}-\frac{1}{d}\sum\limits_{j=1}^{d}f_{x_j}\Big].
\]
Moreover, for all vertices $y\in\gd$ we have
\[
\pi_{z,0}(\delta_y)\ = \  \frac{\ \overline{f_e(y)}\ }{\ \al_0^{2}(z)\ } \,\cdot\, f_e.
\]
\end{f}

\bigskip

\begin{dowod}
Since $\pi_{z,x}(\delta_y)\in A_x$, there exist constants $a_j$ such that
\[
\pi_{z,x}(\delta_y)\ = \ \sum\limits_{j=1}^{d} a_j f_{x_j}.
\]
Let $g=\sum\limits_{j=1}^d b_j f_{x_j}$ be any element of $A_x$.
Then
\[
\iloczyn{g}{\de_y}\ = \ \iloczyn{\pi_{z,x}(g)}{\de_y}\ = \ \iloczyn{g}{\pi_{z,x}(\de_y)}.
\]
Since $\supp (f_{x_j})\subseteq \Gamma_{x_j}$  and $y\in\Gamma_{x_i}$, we obtain
\[
\iloczyn{g}{\de_y} \ =\ \sum\limits_{j=1}^d b_j \,f_{x_j}(y)\ = \ b_i\, f_{x_i}(y).
\]
On the other hand, $\iloczyn{f_{x_i}}{f_{x_j}}=0$ for $i\neq j$. Hence
\[
\iloczyn{g}{\pi_{z,x}(\de_y)}\ = \ \sum\limits_{j=1}^d b_j\, \overline{a_j}\, \|f_{x_j}\|^2\
= \ \al_k^2(z) \cdot \sum\limits_{j=1}^d b_j \, \overline{a_j}.
\]
We thus get the following equation
\[
b_i \, f_{x_i} (y)\ = \ \al_k^2(z)\cdot \sum\limits_{j=1}^d b_j \, \overline{a_j}
\]
for any coefficients $b_j$ such that $\sum\limits_{j=1}^d b_j=0$.
Let $b_i=1$. Then setting $b_{j_o}=-1$ for an arbitrary
$j_0\neq i$ yields
\[
f_{x_i}(y)\ =\ \al_k^2(z)\, (\overline{a_i} -\overline{a_{j_o}}).
\]
It means that  the coefficients  $a_j$ for all $j\neq i$ have the same value as $j_0$ was chosen arbitrarily.
Set $a=a_j$ for $j\neq i$.
Since the sum of all the coefficients $a_j$ vanishes, we have $a_i=-(d-1)a$.
It follows
\[
f_{x_i}(y) \ = \ \al_k^2(z)\cdot \big[  -(d-1)\overline{a} - \overline{a}\, \big]\ =\ -d\,\overline{a}\cdot \al_k^2(z),
\]
whence
\[
a\ =\ -\  \frac{\ \ \overline{f_{x_i}(y)}\ \ }{\ \,d\,\al_k^2(z)\ }.
\]
Summarizing,
\begin{align*}
\pi_{z,x}(\de_y)\
&= \ \frac{(d-1)\, \overline{f_{x_i}(y)} }{ d\, \al_k^2(z) } \cdot f_{x_i}
\ - \
\sum\limits_{j\neq i} \frac{\overline{f_{x_i}(y)}}{\ \,d\,\al_k^2(z)\ } \cdot f_{x_j}\\
&=\ \frac{\ \,\overline{f_{x_i}(y)}\ }{\ \al_k^2(z)\ }
\cdot \bigg[ \Big(1-\frac{1}{d}\Big) \cdot f_{x_i} \ - \ \frac{1}{d}\cdot\sum\limits_{j\neq i} f_{x_j} \bigg]\\
&=\ \frac{\ \,\overline{f_{x_i}(y)}\ }{\ \al_k^2(z)\ }
\cdot \bigg[ \ f_{x_i} \ -\ \frac{1}{d}\cdot \sum\limits_{j=1}^d f_{x_j}\ \bigg].
\end{align*}
The formula for  $\pi_{z,0}(\de_y)$ is clear as
the vector generating the one dimensional subspace $A_0$ is equal to $(\al_0(z))^{-1}f_e$.
Hence
\[
\pi_{z,0}(\de_y)\ = \ \frac{\iloczyn{\de_y}{f_e}}{\|f_e\|^2} \cdot f_e \ =\ \frac{\ \overline{f_e(y)}\ }{\ \al_0^{2}(z)\ } \,\cdot\, f_e.
\]
\end{dowod}

\bigskip
\bigskip
\noindent
In order to describe the  action of $\pi_z$ on $\de_y$ it is necessary to consider
all the subspaces  $A_x$ such that $y\in\Gamma_x$.

\bigskip

\begin{w}\label{pi_z_delta_y}
Assume that the length of the vertex $y\in\gd$ is equal to $n\ge 0$.
Let $y_0, y_1, y_2, \ldots , y_n$ be the path from the root $e$ to $y=y_n$.
Then
\[
\pi_z(\de_y)\ = \  \frac{\,\ \overline{f_e(y)} \ }{\,\  \al_0^2(z) \ } \cdot f_e \ \
+ \ \ \sum\limits_{i=1}^n \frac{\,\ \overline{f_{y_i}(y)} \ }{\,\  \al_i^2(z) \ }
\cdot \Bigg[\ f_{y_i} - \frac{1}{d}\sum\limits_{j=1}^d f_{(y_{i-1})_j}   \ \Bigg],
\]
where $(y_{i-1})_j$ for $j=1,2,\ldots$ are all the successors of the vertex $y_{i-1}$.
\end{w}

\bigskip

\begin{dowod}
It is sufficient to apply Fact \ref{rzuty} to the sum
\[
\pi_z (\de_y)\ = \ \sum\limits_{x\in\gd\cup\{0\}}  \pi_{z,x} \, (\de_y)\ = \
\big(\,\pi_{z,y_0}+\pi_{z,y_1}+\ldots+\pi_{z,y_n} \big)\, (\de_y).
\]

\end{dowod}

\bigskip
\bigskip

\noindent
Note that in view of (\ref{odpowiednie}) and (\ref{odpowiednie2}),
the isometry $U$ defined by (\ref{izometria}) can be expressed  by the following formula
\[
U \bigg(\ \al_{k}(z)\, \sqrt{d^{k}}\cdot \sum\limits_{i=1}^{d} a_i \cdot \1{\Omega_{x_i}}\,\bigg) \ = \ \
 \ \sum\limits_{i=1}^d a_i \, f_{x_i},
\]
where $k=|x_i|=|x|+1$, and
\begin{equation}\label{U_0}
U \big(\, \al_0(z)\cdot \1{\Omega} \,\big)\  = \  f_e.
\end{equation}
As the supports are disjoint for $i\neq j$, the first formula can be written down as follows
\begin{equation}\label{U}
U \bigg(\ \al_{k}(z)\, \sqrt{d^{k}}\cdot  \1{\Omega_{x_i}}\,\bigg) \ = \ \
 \  f_{x_i}.
\end{equation}

\bigskip
\bigskip
Here comes the promised explicit formula for the Poisson kernel.
\bigskip
\bigskip

\begin{tw}
Let  $\omega\in\Omega$ and   $y\in\gd$ be of length $n$. Then
\[
P_z(y,\omega) \ = \
\frac{\ f_e(y)\ }{\ \al_0(z)\ } \cdot\1{\Omega} \ + \
\sum\limits_{i=1}^n
\frac{\,\ {f_{y_i}(y)} \ }{\,\  \al_i(z) \ }\cdot \sqrt{d^i} \cdot \Big[\ \1{{\Omega_{y_i}}} - \frac{1}{d} \1{\Omega_{y_{i-1}}}   \ \Big],
\]
where $\{ y_0, y_1, y_2, \ldots , y_n\}$ is the only one path from the root $e=y_0$ to the vertex $y=y_n$.
\end{tw}

\bigskip

\begin{dowod}
Applying (\ref{U_0}) and (\ref{U}) to Corollary \ref{pi_z_delta_y} yields
\[
\pi_z(\de_y)\ = \
U \big( \ S(y)\ \big),
\]
where
\[
S(y) \ = \ \frac{\ \overline{f_e(y)}\ }{\ \al_0(z)\ } \cdot\1{\Omega} \ + \
\sum\limits_{i=1}^n
\frac{\,\ \overline{{f_{y_i}(y)}}\sqrt{d^i} \ }{\,\  \al_i(z) \ } \cdot \Big[\ \1{{\Omega_{y_i}}} - \frac{1}{d} \1{\Omega_{y_{i-1}}}   \ \Big].
\]

\noindent
For any function $F\in \Ldwa$ we have
\[
U F (y) \ = \ \iloczyn{ UF}{\de_y }\ = \ \iloczyn{ UF}{\pi_z(\de_y)}.
\]
On the other hand, the Poisson kernel $P_z$ satisfies (\ref{riesz}) so
\[
UF (y)\ = \ \int\limits_{\Omega} P_z(y,\omega) F(\omega) \dw \ = \ \iloczyn{F}{\overline{P_z(y, \cdot)}}.
\]
By the above, we obtain
\begin{align*}
\iloczyn{F}{\overline{P_z(y, \cdot)}} \ &= \ \iloczyn{ UF}{\pi_z(\de_y)} \\ &= \
\iloczyn{UF}{U(S(y))}\
= \ \iloczyn{F}{S(y)},
\end{align*}
which completes the proof.
\end{dowod}



\hspace{7pt}

\begin{center}
{\textsc{The spectrum of a selfadjoint extension}}
\end{center}

\hspace{7pt}

Our next aim is to describe the spectral properties of $J$.

\bigskip

Recall that we are considering the case when $J$ is not essentially selfadjoint.
For a fixed vertex $x\in\gd$ we define the linear subspace $H_x$
of the Hilbert space $\ell^2(\gd)$ consisting of those functions $f\in  \ell^2(\gd)$ which satisfy
\begin{itemize}
    \item[{\scriptsize{${(1)}$}}]  $\ \supp(f)\subset \Gamma_x \setminus\{x\}$,
    \item[{\scriptsize{${(2)}$}}]  $\ \sum\limits_{i=1}^{d} f(x_i)=0$,
    \item[{\scriptsize{${(3)}$}}]  $f$ is radial on each subtree $\Gamma_{x_i}$,
    \item[{\scriptsize{${(4)}$}}]  the value on a level of $\Gamma_{x_i}$ is proportional to the value on the corresponding level of $\Gamma_{x_j}$
    with the coefficient $\frac{f(x_i)}{f(x_j)}$.
\end{itemize}
Moreover, we set $H_0=\ell^2_r(\gd)$.

\bigskip
\bigskip

\begin{f}
The family of the spaces $H_x\subset \ell^2(\gd)$, where $x\in\gd\cup\{0\}$, satisfies
\begin{itemize}
    \item[${(1)}$] $J[H_x]\,\subseteq \, H_x$ for every $x$,
    \item[${(2)}$] $H_x$ is closed for every $x$,
    \item[${(3)}$] $H_x\,\perp\, H_y$ for $x\neq y$,
    \item[${(4)}$] $\bigoplus\limits_{x\in\gd\cup\{ 0\}} H_x\,=\,\lin\{f\in H_x\colon \ x\in\gd\cup\{0\}\} $ is dense in $\ell^2(\gd)$.
\end{itemize}
\end{f}

\bigskip

\begin{dowod}
Properties (1) and (2) are clear. Property (3) may be proved in much the same way as  Fact \ref{ortogonalnosc}.
In order to prove (4) let us assume that $g\in\ell^2(\gd)$ is orthogonal to every $H_x$ for $x\in\gd\cup \{0\}$.
In particular, by the orthogonality to $\de_e\in H_0$, we obtain
\[
0\ = \ \iloczyn{g}{\de_e}\ = \ g(e).
\]
The orthogonality to $\de_{e_i}-\de_{e_j}\in H_e$ gives
\[
0\ = \ \iloczyn{g}{\de_{e_i}-\de_{e_j}}\ = \ g(e_i) \, -\, g(e_j),
\]
whence  the values on the first level are all equal.
Furthermore, since the characteristic function $\chi_1$ of the first level is an element of $H_0$, we get
\[
0\ = \ \iloczyn{g}{\chi_1}\ = \ \sum\limits_{i=1}^d g(e_i).
\]
It means that $g$ vanishes also at the first level of $\gd$.
Similar considerations show that $g$ is equal to 0 at each level of $\gd$.
\end{dowod}

\bigskip
\bigskip

\noindent
Let $J_x$ denote the restriction of $J$ to the subspace $H_x\cap D(J)$.
For $x=0$  the operator $J_x$ is expressed by the matrix $J^r=J^r_0$.
For the vertex $x\in\gd$ the action of $J_x$ is associated with the restricted matrix $J^r_n$, where $n=|x|+1$  (cf. (\ref{macierz_obcieta})).

\bigskip

\noindent
Since $J$ is not essentially selfadjoint, neither is any of the matrices $J^r_n$
(cf. the beginning of part \emph{The Description of the deficiency space}).
It is known that there exists a selfadjoint extension $\tilde{J_n}$ for such a matrix $J^r_n$
and the spectrum of each selfadjoint extension is a discrete set
(cf. Theorem \ref{n_ekstremalne}).

\bigskip

\noindent
Let $\tilde{J_x}$ be the operator with the domain $D(\tilde{J_x})\subseteq H_x$
associated with the selfadjoint extension $\tilde{J}_{|x|+1}$.
Hence its spectrum
$\sigma(\tilde{J}_{|x|+1})$ is a discrete set so $\tilde{J_x}$ has a pure point spectrum
(i.e. there exists a basis consisting of eigenvectors).
Define $\tilde{J}$ by
\[
\tilde{J}(f)\ = \ \sum\limits_{x\in\gd\cup\{0\}}\tilde{J}_x(f)
\]
with domain
\[
D(\tilde{J})\ = \ \bigoplus D(\tilde{J}_x)\ = \ \lin \{ f\in\ldwa\colon \ f\in D(\tilde{J}_x)  \,\}.
\]
Since the $H_x$ are invariant under $J$ and the Hilbert orthogonal sum of them is equal to the whole
space $\ell^2(\gd)$, the operator $\tilde{J}$ is a selfadjoint extension of $J$.
Moreover, the spectrum of this extension
\[
\sigma(\tilde{J})\ =\ \overline{\,\bigcup\, \sigma\big(\tilde{J}_x\big)}
\]
is also a pure point spectrum.


\clearpage

\section{\textsc{A JACOBI OPERATOR ON A TREE WITH ONE END}}

\hspace{7pt}

For a fixed number $d=2,3,4,\ldots$  we consider an infinite homogeneous tree
of degree $d$ which is partially ordered and locally looks like the one in the former chapter but upside down.
For instance, if $d=3$, the top levels of the tree look as follows

\bigskip

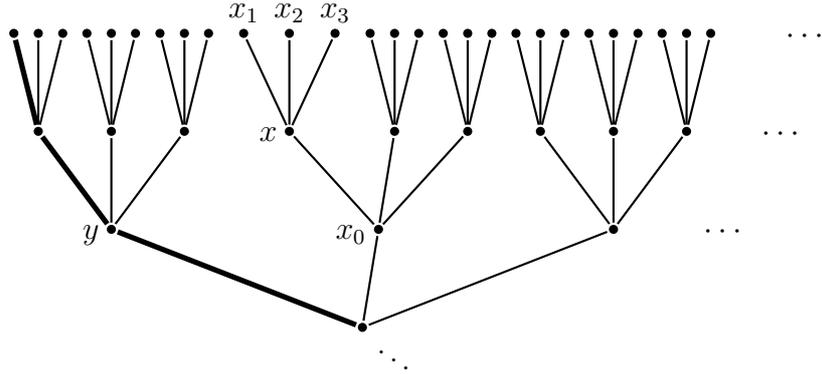
\begin{center}
\newcommand{\dedge}{\ncline[%
linewidth=2pt%
]}
\pstree[nodesep=0.3mm,treemode=U,treefit=loose,treesep=0.2,levelsep=1.3cm]
{\Tc*{0.6mm}~[tnpos=b,tnsep=0mm]{$\qquad\ddots$}}
{

\pstree
{\Tc*[edge=\dedge]{0.6mm}~[tnpos=l,tnsep=1mm]{$y  $}}
{
  \pstree
  {\Tc*[edge=\dedge]{0.6mm}~[tnpos=l,tnsep=1mm]{$ $}}
  {
   \Tc*[edge=\dedge]{0.6mm}~[tnpos=b,tnsep=1mm]{$ $}
   \Tc*{0.6mm}~[tnpos=b,tnsep=1mm]{$ $}
   \Tc*{0.6mm}~[tnpos=b,tnsep=1mm]{$ $}
  }
    \pstree
     {\Tc*{0.6mm}~[tnpos=l,tnsep=1mm]{$ $}}
     {  \Tc*{0.6mm}~[tnpos=b,tnsep=1mm]{$ $}
        \Tc*{0.6mm}~[tnpos=b,tnsep=1mm]{$ $}
        \Tc*{0.6mm}~[tnpos=b,tnsep=1mm]{$ $}
     }
       \pstree
        {\Tc*{0.6mm}~[tnpos=l,tnsep=1mm]{$ $}}
        {
        \Tc*{0.6mm}~[tnpos=b,tnsep=1mm]{$ $}
        \Tc*{0.6mm}~[tnpos=b,tnsep=1mm]{$ $}
        \Tc*{0.6mm}~[tnpos=b,tnsep=1mm]{$ $}
  }
}

 \pstree
{\Tc*{0.6mm}~[tnpos=l,tnsep=1mm]{$ x_0 $}}
{
  \pstree
  {\Tc*{0.6mm}~[tnpos=l,tnsep=1mm]{$x $}}
  {
   \Tc*{0.6mm}~[tnpos=a,tnsep=0mm]{$x_1$}
   \Tc*{0.6mm}~[tnpos=a,tnsep=0mm]{$x_2 $}
   \Tc*{0.6mm}~[tnpos=a,tnsep=0mm]{$x_3 $}
  }
    \pstree
     {\Tc*{0.6mm}~[tnpos=l,tnsep=1mm]{$ $}}
     {  \Tc*{0.6mm}~[tnpos=b,tnsep=1mm]{$ $}
        \Tc*{0.6mm}~[tnpos=b,tnsep=1mm]{$ $}
        \Tc*{0.6mm}~[tnpos=b,tnsep=1mm]{$ $}
     }
       \pstree
        {\Tc*{0.6mm}~[tnpos=l,tnsep=1mm]{$ $}}
        {
        \Tc*{0.6mm}~[tnpos=b,tnsep=1mm]{$ $}
        \Tc*{0.6mm}~[tnpos=b,tnsep=1mm]{$ $}
        \Tc*{0.6mm}~[tnpos=b,tnsep=1mm]{$ $}
  }
}

\pstree
{\Tc*{0.6mm}~[tnpos=r,tnsep=3mm]{$\qquad\cdots $}}
{
  \pstree
  {\Tc*{0.6mm}~[tnpos=r,tnsep=1mm]{$ $}}
  {
   \Tc*{0.6mm}~[tnpos=b,tnsep=1mm]{$ $}
   \Tc*{0.6mm}~[tnpos=b,tnsep=1mm]{$ $}
   \Tc*{0.6mm}~[tnpos=b,tnsep=1mm]{$ $}
  }
    \pstree
     {\Tc*{0.6mm}~[tnpos=l,tnsep=1mm]{$ $}}
     {  \Tc*{0.6mm}~[tnpos=b,tnsep=1mm]{$ $}
        \Tc*{0.6mm}~[tnpos=b,tnsep=1mm]{$ $}
        \Tc*{0.6mm}~[tnpos=b,tnsep=1mm]{$ $}
     }
       \pstree
        {\Tc*{0.6mm}~[tnpos=r,tnsep=1mm]{$\qquad\cdots $}}
        {
        \Tc*{0.6mm}~[tnpos=b,tnsep=1mm]{$ $}
        \Tc*{0.6mm}~[tnpos=b,tnsep=1mm]{$ $}
        \Tc*{0.6mm}~[tnpos=r,tnsep=1mm]{$\qquad\cdots$}
  }
}
}.
\end{center}

\bigskip

In view of the picture above, it is intuitively clear what the partial order  in this tree is.
All vertices with only one edge are on the zero level.
Those vertices which are at a distance 1 from the zero level, have length 1.
And so on.
To be more precise, this time we distinguish not a vertex but an infinite path
$\omega=\{\omega_0,\omega_1,\omega_2,\ldots\}$ where $\omega_0$ is any vertex with only one edge.
The natural distance $\dist(\cdot,\cdot)$ enables one to calculate
the distance between a given vertex $x$ and the path $\omega$, i.e. $\dist(x,\omega)$.
Then, by the {\emph{length}} of a vertex $x$ we mean the following difference
\[
|x|\ =\ n - \dist(x,\omega),
\]
where $n$ is equal to the index of the element of $\omega$ which realizes the distance
\[
\dist(x,\omega)\ = \ \dist(x,\omega_n).
\]
In the picture above the fixed path $\omega$ is indicated by a bold line.
For the marked vertex $x$ we have $|x|=3-2=1$ and for $y=\omega_2$ we have $|y|=2-0=2$.
It is clear that the length $| \cdot |$ defined in this way is independent  of the choice of $\omega$.

\bigskip

The set of all vertices with the defined partial order is denoted by $\Lambda_d$.

\bigskip

In the tree $\Lambda_d$ each vertex of length at least equal to  1 has exactly $d$ predecessors and 1 successor.
Each origin, i.e. the vertex with no predecessor, has length 0 and exactly one edge (downward) so also 1 successor.
This time there are infinitely many vertices of length 0.
At each vertex, however, there is just one edge downward so $\Lambda_d$ can be said to
be a homogeneous tree with one end.

\bigskip

Analogously to $\gd$, the predecessors of a  vertex $x$ are denoted by $x_1,x_2,\ldots,x_d$ and the successor by $x_0$.
In analogy with the previous chapter we also define the action of the Jacobi operator $J$
on the characteristic function $\delta_x$  of a vertex $|x|=n$, namely
\[
J\delta_x \ = \ \la_{n-1}\big(\delta_{x_1}+\delta_{x_2}+\ldots+\delta_{x_d}\big)+\beta_n\delta_x+\lambda_n\delta_{x_0}.
\]
The domain of $J$ consists of functions with finite supports, i.e.
\[
D(J)\ =\ \lin\big\{\delta_x\colon\quad x\in\Lambda_d \big\}\ \subseteq\ \ell^2(\Lambda_d).
\]
We still keep the convention that $\lambda_{-1}=0$ what makes the formula for $J$ clear also for the vertices of length 0.

\bigskip
\bigskip

\begin{f}\label{jadro}
The deficiency space $N_z(J)$ of the operator $J$ on $\ell^2(\ld)$ consists of all square-summable functions on $\ld$
satisfying
\begin{equation}\label{rec}
z v(x)\,=\,\lambda_{n-1}\big(v({{x}_1}) + \ldots +v({{x}_d})\big)+ \beta_nv(x) + \lambda_n v({x}_0)
\end{equation}
for all $|x|=n$ and  all $n\geq 0$.
\end{f}

\begin{dowod}
In analogy with  $\gd$  the assertion is implied by the following calculation:
\begin{eqnarray*}
0 &=&(v\;,\; (J-\bar{z}) \de_x)  \\
  &=&(v\;,\; \lambda_{n-1}\big(\de_{{x}_1} + \ldots +\de_{{x}_d}\big) + \beta_n\de_x
               + \lambda_n \de_{{x}_o}
               -\bar{z}\de_x)  \\
  &=& \lambda_{n-1}\big(v({{x}_1}) + \ldots +v({{x}_d})\big)+ \beta_nv(x)
               + \lambda_n v({x}_0)
               -z v(x).
\end{eqnarray*}
\vspace{6pt}
\end{dowod}

\noindent
{\textbf{Remark.}} Clearly, an equivalent formulation of the assertion is:
\[
N_z(J)\ =\ \Big\{\, v\in\ell^2(\ld)\colon \ \ \ Jv\, (x)\ =\ z \cdot v(x), \qquad x\in\Lambda_d\, \Big\}.
\]

\bigskip

\noindent
Let $\Lambda_x$ denote the subtree of $\ld$ which ends at the vertex $x$. The subtree $\Lambda_x$ is marked in the picture below.

\bigskip

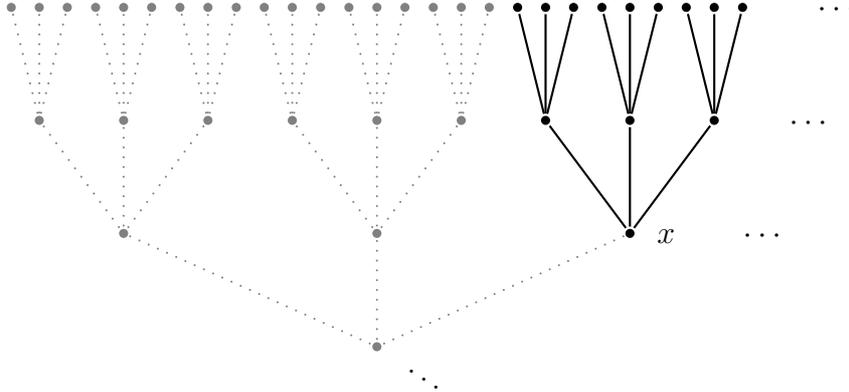
\begin{center}
\pstree[nodesep=0.3mm,treemode=U,treefit=loose,treesep=0.25,levelsep=1.5cm,linestyle=dotted,linecolor=gray]
{\Tc*[linecolor=gray]{0.6mm}~[tnpos=b,tnsep=0mm]{$\qquad\quad\ddots$}}
{

\pstree
{\Tc*{0.6mm}~[tnpos=l,tnsep=3mm]{$  $}}
{
  \pstree
  {\Tc*{0.6mm}~[tnpos=l,tnsep=1mm]{$ $}}
  {
   \Tc*{0.6mm}~[tnpos=b,tnsep=1mm]{$ $}
   \Tc*{0.6mm}~[tnpos=b,tnsep=1mm]{$ $}
   \Tc*{0.6mm}~[tnpos=b,tnsep=1mm]{$ $}
  }
    \pstree
     {\Tc*{0.6mm}~[tnpos=l,tnsep=1mm]{$ $}}
     {  \Tc*{0.6mm}~[tnpos=b,tnsep=1mm]{$ $}
        \Tc*{0.6mm}~[tnpos=b,tnsep=1mm]{$ $}
        \Tc*{0.6mm}~[tnpos=b,tnsep=1mm]{$ $}
     }
       \pstree
        {\Tc*{0.6mm}~[tnpos=l,tnsep=1mm]{$ $}}
        {
        \Tc*{0.6mm}~[tnpos=b,tnsep=1mm]{$ $}
        \Tc*{0.6mm}~[tnpos=b,tnsep=1mm]{$ $}
        \Tc*{0.6mm}~[tnpos=b,tnsep=1mm]{$ $}
  }
}

 \pstree
{\Tc*{0.6mm}~[tnpos=l,tnsep=3mm]{$  $}}
{
  \pstree
  {\Tc*{0.6mm}~[tnpos=l,tnsep=1mm]{$ $}}
  {
   \Tc*{0.6mm}~[tnpos=b,tnsep=1mm]{$ $}
   \Tc*{0.6mm}~[tnpos=b,tnsep=1mm]{$ $}
   \Tc*{0.6mm}~[tnpos=b,tnsep=1mm]{$ $}
  }
    \pstree
     {\Tc*{0.6mm}~[tnpos=l,tnsep=1mm]{$ $}}
     {  \Tc*{0.6mm}~[tnpos=b,tnsep=1mm]{$ $}
        \Tc*{0.6mm}~[tnpos=b,tnsep=1mm]{$ $}
        \Tc*{0.6mm}~[tnpos=b,tnsep=1mm]{$ $}
     }
       \pstree
        {\Tc*{0.6mm}~[tnpos=l,tnsep=1mm]{$ $}}
        {
        \Tc*{0.6mm}~[tnpos=b,tnsep=1mm]{$ $}
        \Tc*{0.6mm}~[tnpos=b,tnsep=1mm]{$ $}
        \Tc*{0.6mm}~[tnpos=b,tnsep=1mm]{$ $}
  }
}

\pstree[linecolor=black,linestyle=solid]
{\Tc*[linecolor=black]{0.6mm}~[tnpos=r,tnsep=3mm]{$ x \qquad\cdots $}}
{
  \pstree
  {\Tc*{0.6mm}~[tnpos=l,tnsep=1mm]{$ $}}
  {
   \Tc*{0.6mm}~[tnpos=b,tnsep=1mm]{$ $}
   \Tc*{0.6mm}~[tnpos=b,tnsep=1mm]{$ $}
   \Tc*{0.6mm}~[tnpos=b,tnsep=1mm]{$ $}
  }
    \pstree
     {\Tc*{0.6mm}~[tnpos=l,tnsep=1mm]{$ $}}
     {  \Tc*{0.6mm}~[tnpos=b,tnsep=1mm]{$ $}
        \Tc*{0.6mm}~[tnpos=b,tnsep=1mm]{$ $}
        \Tc*{0.6mm}~[tnpos=b,tnsep=1mm]{$ $}
     }
       \pstree
        {\Tc*{0.6mm}~[tnpos=r,tnsep=1mm]{$ \qquad\cdots $}}
        {
        \Tc*{0.6mm}~[tnpos=b,tnsep=1mm]{$ $}
        \Tc*{0.6mm}~[tnpos=b,tnsep=1mm]{$ $}
        \Tc*{0.6mm}~[tnpos=r,tnsep=1mm]{$\qquad\cdots$}
  }
}

}
\end{center}

\bigskip
\bigskip
The following technical lemma is a direct preparation for the main theorem
of this chapter which considers essential selfadjointness of~$J$ and will follow next.
\bigskip
\bigskip

\begin{lem}{\label{techniczny}}
Assume that for a number $z\in\C$ (we allow $z$ to be real) a function $v\in\ltwo$ satisfies
the recurrence relation (\ref{rec}) describing $N_z(J)$.
Let a vertex $x\in\ld$ have length $n$.
Then the values of $v$ on $\Lambda_x$ are constant on each level of this subtree.
Moreover, if $y\in\Lambda_x$ and if $|y|=k\ge 0$, then
\[
v(y) \ = \ \sqrt{d^k}\,p_k(z) \cdot v_0,
\]
where $v_0$ is the value of $v$ on the zero level of $\Lambda_x$ and the numbers
$p_n(z)$ are the values of the orthogonal polynomials associated with the matrix $J^r$ (cf. (\ref{p_n})).
\end{lem}

\bigskip

\begin{dowod}
The proof is by induction on $n$. \\
(1) Fix a vertex $|x|=n=1$. For each $i=1,2,3,\ldots ,d$ we have by (\ref{rec}) for $x=x_i$
\begin{align*}
(z-\beta_0) v(x_i)\ = \ \lambda_0 v(x).
\end{align*}
Hence
\[
v(x)\ =\ \frac{z-\beta_0}{\la_0} \cdot v(x_i) \ = \ \sqrt{d}\,p_1(z) \cdot v_0
\]
as $p_0(z)=1$ and
\[
z\, p_0(z)\ = \ \beta_0 p_0(z)\ +\ \sqrt{d}\,\la_0\, p_1(z).
\]
(2) Assume that the assertion holds for some $n\ge 1$.
Let $x$ be any vertex in $\Lambda_x$ of length $n+1$.
Each of its predecessors $x_i$ has length $n$ so the values of $v$ on the $k$th level
of $\Lambda_{x_i}$ are constant and equal to $\sqrt{d^k} \cdot p_k(z)\cdot v_0^i$ respectively,
where $v_0^i$ is the value of $v$ on the zero level of~$\Lambda_{x_i}$.
By assumption, the recurrence equation (\ref{rec}) at~$x_i$
\[
\big(z - \beta_{n}\big)\cdot v(x_i)\ =\ \la_{n-1}\cdot\sum\limits_{j=1}^d v\big({(x_i)}_j\big) \,+\,\la_{n} \cdot v(x)
\]
yields
\[
\big(z - \beta_{n}\big)\cdot \sqrt{d^n}\, p_n(z) \cdot v_0^i\ =\ d\la_{n-1}\cdot \sqrt{d^{n-1}}\, p_{n-1}(z) \cdot v_0^i \,+\,\la_{n} \cdot v(x).
\]
Hence
\[
v(x) \ = \  \frac{(z - \beta_{n}) \sqrt{d^n}\,p_n(z)-d\la_{n-1} \sqrt{d^{n-1}}\,p_{n-1}(z) }{\la_n} \cdot v_0^i.
\]
By the recurrence relation (\ref{gwiazdka}) satisfied by $\{p_n(z) \}$ we get
\[
v(x)\ = \ v_0^i \cdot \sqrt{d^{n+1}} \cdot p_{n+1}(z).
\]
\end{dowod}

\bigskip
\bigskip
Here is the main theorem.
\bigskip
\bigskip

\begin{tw}\label{main_2}
The operator $J$ on $\ld$ is always essentially selfadjoint.
\end{tw}

\bigskip

\begin{dowod}
For a complex number $z\notin\R$ all the coefficients appearing in Lemma \ref{techniczny}, i.e. the numbers
\[
\sqrt{d^k}\,p_k(z),
\]
are nonzero as all roots of the orthogonal polynomials $p_n$ are real.
By Lemma \ref{techniczny}, if there existed  a function satisfying all the recurrence equations (\ref{rec})
describing $N_z(J)$,
it would have to be nonzero and constant on levels of the whole tree $\ld$.
However, there are infinitely many vertices on each level
so such a function cannot be square-summable.
Therefore,  $N_z(J)=\{ 0\}$.
\end{dowod}

\bigskip

\begin{tw}
The Jacobi operator $J$ on $\ld$ has a pure point spectrum $\sigma (J)$, i.e. there is an orthonormal basis
consisting of eigenvectors for $J$.
Moreover, $\sigma (J)$ coincides with the closure of
the set of all roots of the orthogonal polynomials $p_n$ associated with the matrix $J$.
\end{tw}

\bigskip

\begin{dowod}
Since $J$ is essentially selfadjoint it suffices to point at a set of eigenvectors which is linearly dense in $D(J)$.

Fix a vertex $x\in\ld$ of length $n\ge 1$.
We consider a subspace $M_x \subset D(J)$ consisting of those functions whose supports
are contained in $\Lambda_x$.
Clearly,
\[
\dim M_x \ = \ 1+ d+d^2+\ldots+d^n.  
\]
It is known that the polynomial $p_{n}$ has exactly $n$ real simple roots
\[
t_1,\ t_2,\ t_3,\ \ldots ,\ t_{n}.
\]
For a fixed predecessor  $x_i$ of $x$ and  for a fixed root  $t_j$ of $p_n$
let $f_{i,j}\in M_x$ be given by
\[
f_{i,j}(y)\ = \ \left\{
\begin{array}{ll}
\sqrt{d^k}\cdot p_k(t_j)\ \quad & \textrm{ for }\quad y\in\Lambda_{x_i}\ \textrm{ and }\ |y|=k,\\
0 & \textrm{ for } y\notin \Lambda_{x_i}.
\end{array}\right.
\]
Of course, thus defined $f_{i,j} $ satisfies the recurrence equations (\ref{rec}) (note that it is for $z=t_j\in\R$)
contained in the subtree $\Lambda_{x_i}$.
Furthermore, since
\[
f_{i,j}(x)\ = \ 0 \ = \ \sqrt{d^{n}} \cdot p_{n}(t_j),
\]
the recurrence equation (\ref{rec}) at $x_i$ holds also.
Hence, the linear combinations
\[
f_{1,j}\ -\ f_{i,j} \qquad \textrm{ for }\qquad i=2,3,\ldots,d
\]
satisfy, in addition, the recurrence equation (\ref{rec}) at $x$, i.e.
\[
0=(z-\beta_n)(f_{1,j}(x)-f_{i,j}(x))= \la_{n-1}(f_{1,j}(x_1)-f_{i,j}(x_i))+\la_n\cdot 0
\]
because $f_{1,j}(x_1)=f_{i,j}(x_i)$.

\bigskip

By the above, when $j$ is fixed and $i$ varies from 2 to $d$ the functions $f_{1,j}-f_{i,j}$ satisfy
the recurrence equations (\ref{rec}) for $z=t_j$ at every vertex  of the tree $\ld$, i.e.
\[
J\,\big(f_{1,j}-f_{i,j}\big)\ = t_j \cdot \big(f_{1,j}-f_{i,j}\big),\qquad i=2,3,\ldots,d.
\]
Hence they are  eigenfunctions associated with the eigenvalue $t_j$.
Clearly, there are $d-1$ of them and they form a linearly independent system since
functions $f_{i,j}$ are pairwise orthogonal for $i=1,2,\ldots,d$ as functions with disjoint supports.

\bigskip

In this way, for a fixed vertex of length $n$,  we indicated exactly $n\cdot (d-1)$ linearly independent eigenfunctions
associated with this vertex.
If we consider the entire subtree $\Lambda_x$  there are  $d^{n-k}$ vertices of given length $k$.
Of course, the eigenfunctions corresponding to two such vertices of given length $k$ are orthogonal
as their supports are disjoint.
Moreover, considering two vertices such that one is in the subtree associated with another one,
the corresponding eigenfunctions are also orthogonal. This is because on each level
of the smaller tree one function has a constant value while the values of the other one
sum up to zero (cf. the proof of Fact \ref{ortogonalnosc}).
Therefore, the number of all thus defined eigenfunctions for $J$ with the supports contained in the subtree $\Lambda_x$ is equal to
\[
(d-1)\cdot\sum\limits_{k=1}^{n} k\cdot d^{n-k}
\]
and all of them form a linearly independent system.

\bigskip

Let $V_x \subseteq M_x$ denote the linear subspace  spanned by the  eigenvectors defined above and with support contained in $\Lambda_x$.
Then
\[
\dim V_x =  (d-1)\cdot\sum\limits_{k=1}^{n} k\cdot d^{n-k}=\big( 1+d+d^2+\ldots +d^n \big) -(n+1)  .
\]

Since there are $n+1$ levels in $\Lambda_x$, there exists exactly $n+1$ linearly independent functions in $M_x$
which are constant on the levels of $\Lambda_x$.
Therefore, the equality
\[
\dim M_x \ = \ \dim V_x\  + \ (n+1),
\]
obtained above, means that the orthogonal complement of $V_x$ in $M_x$ consists only of functions constant on levels of
$\Lambda_x$.

\bigskip

To complete the proof it suffices to show that there is no  square-summable and nonzero functions which are orthogonal to every $V_x$.
Assume that $f\in\ltwo$ satisfies
\[
\forall x\in\Lambda_d \qquad  f\perp V_x.
\]
Then $f$ is constant on levels of $\Lambda_x$ for each vertex $x\in\Lambda_d$.
Hence $f$ is constant on levels of the whole tree $\Lambda_d$.
But $f$ is square-summable. Therefore $f\equiv 0$.
\end{dowod}


\clearpage
\thispagestyle{empty}


\begin{thebibliography}{9}
\bigskip
\normalsize

\bibitem{akhiezer} N. I. Akhiezer, \emph{The Classical Moment Problem}, Hafner Publishing Co., New York, 1965.

\bibitem{berg} C. Berg, \emph{Moment problems and orthogonal polynomials}, lecture notes, 2007.

\bibitem{chihara} T. Chihara, \emph{An Introduction to Orthogonal Polynomials}, Mathematics and Its Applications,
Vol. 13, Gordon and Breach, New York, London, Paris, 1978.

\bibitem{hamburger} H. Hamburger, \emph{\"Uber Eine Erweiterung der Stieltjesschen Momentenproblems},
Math. Ann. 81, 235--319; 82, 120--164; 82, 168--187, 1920, 1921.

\bibitem{nevanlinna} R Nevanlinna, \emph{Asymptotische Entwicklungen beschr\"ankter Funktionen und das Stieltjessche Momentenproblem},
Ann. Acad. Scient. Fennicae (A) 18, no. 5 (52 pp.), 1922.

\bibitem{reedsimon1} M. Reed, B. Simon, \emph{Methods of Modern Mathematical Physics, I. Functional Analysis}, Academic Press,
New York, 1972.

\bibitem{reedsimon2} M. Reed, B. Simon, \emph{Methods of Modern Mathematical Physics, II. Fourier Analysis, Self-Adjointness}, Academic Press,
New York, 1975.

\bibitem{riesz1} M. Riesz, \emph{Sur le probl\`eme des moments. Premi\`ere Note}, Arkiv f\"or matematik, astronomi och fysik 16, no. 12 (23 pp.), 1921.

\bibitem{riesz2} M. Riesz, \emph{Sur le probl\`eme des moments. Deuxi\`eme Note}, Arkiv f\"or matematik, astronomi och fysik 16, no. 19 (21 pp.), 1922.

\bibitem{riesz3} M. Riesz, \emph{Sur le probl\`eme des moments. Troisi\`eme Note}, Arkiv f\"or matematik, astronomi och fysik 17, no. 16 (62 pp.), 1923.

\bibitem{simon} B. Simon, \emph{The classical moment problem as self-adjoint finite difference operator}, Advances in Mathematics 137 (1998), 82-203.

\bibitem{stieltjes} T. Stieltjes, \emph{Recherches sur les fractions continues}, Ann. Ec. Sci. Univ. Toulouse 8 J1--J122; 9 A5--A47, 1894.

\bibitem{stone} M. Stone, \emph{Linear transformations in Hilbert space and their applications to analysis}, colloq. Publ. Vol. 15,
Amer. Math. Soc., New York, 1932.

\bibitem{szwarc} R. Szwarc, \emph{Wielomiany ortogonalne i problem
moment\'ow}, lecture notes, University of Wroc\l aw, Institute of Mathematics, 2004, http://www.math.uni.wroc.pl/~szwarc/pdf/momentw.pdf
\end{thebibliography}
\end{document}